\documentclass{imsart}

\RequirePackage[OT1]{fontenc}
\RequirePackage{amsthm,amsmath,amssymb}
\usepackage{dsfont}
\RequirePackage[numbers]{natbib}
\RequirePackage[colorlinks,citecolor=blue,urlcolor=blue]{hyperref}
\usepackage{graphicx}

\startlocaldefs

\newcommand \id{\mathds 1}
\newcommand{\Cov}{\mathrm{Cov}}
\newcommand {\R} {\mathbb{R}}
\newcommand {\C} {\mathbb{C}}
\newcommand {\Z} {\mathbb{Z}}
\newcommand{\E}{\mathbb{E}}
\renewcommand {\P} {\mathbb{P}}
\newcommand{\F}{\mathcal{F}}

\mathchardef\mhyphen="2D
\newcommand{\brb}[1]{\left[#1\right]}

\newtheorem{theorem}{Theorem}[section]

\newtheorem{question}[theorem]{Question}

\newtheorem{conjecture}[theorem]{Conjecture}
\newtheorem{example}[theorem]{Example}

\theoremstyle{remark}
\newtheorem{remark}[theorem]{Remark}

\numberwithin{equation}{section}

\endlocaldefs

\begin{document}

\begin{frontmatter}

\title{Gaussian fields and percolation}
\runtitle{Gaussian fields and percolation}


\begin{aug}
\author{\fnms{Dmitry} \snm{Beliaev}\ead[label=e1]{belyaev@maths.ox.ac.uk}}


\runauthor{D. Beliaev}

\end{aug}



\begin{abstract}
In the last two decades there was a lot of progress in understanding the geometry of smooth Gaussian fields. This survey aims to cover one particular line of research: the large scale behaviour of level and excursion sets and their (conjectured) connection to the percolation theory.
\end{abstract}





\end{frontmatter}


\section{Introduction}

Gaussian fields appear naturally in many areas of science (e.g. in quantum chaos \cite{JaSa}, medical imaging \cite{Worsley}, oceanography \cite{CoxMunk,Longuet-Higgins}, cosmology \cite{Bardeen} etc.). They are also very interesting objects from a purely mathematical point of view.

Our goal is to survey recent results aiming at understanding the large scale structure of level sets of Gaussian fields and how it is connected to the percolation theory.

\subsection{Gaussian fields}
We start with a very brief introduction to Gaussian fields, mostly in the stationary case. We do not aim at comprehensive coverage. Moreover, most of this section is not even completely rigorous it should be treated more like a heuristics that is ultimately correct but requires some extra work.  Our goal here is just to explain how one could and should think of these fields. We refer to proper books \cite{AdTa,AzWsch,Bogachev,Lifshits} for more information about Gaussian fields and their basic properties.

There are different ways to think about Gaussian fields. The first one is `probabilistic'. One says that a Gaussian field on an index set $X$ is a collection of jointly Gaussian random variables $F_x=F(x)$ indexed by elements of $X$. It is a well known fact that such collections can be uniquely determined by their mean values $m(x)=\E[F(x)]$ and pairwise covariances $K(x,y)=\E[F(x)F(y)]$. Since the mean value can be changed to $0$ by subtracting the deterministic function $m(x)$, for the rest of this paper we assume that all fields are centred i.e.  $m(x)=0$ for all $x$.  In this definition of the field the nature of the index set $X$ is completely irrelevant, but it becomes important when we ask whether such fields exist and whether they could be constructed from  $K(x,y)$. Kolmogorov's  heorem gives the existence of such a field. It can be shown that if $K$ is sufficiently `nice' then the field is also `nice'. In particular, the modulus of continuity of $F$ can be explicitly written in terms of $K$. We refer to \cite[Section 1.3]{AdTa} for precise statements and complete proofs. Another exposition that we really like is the Appendix in \cite{NaSo_asymptotic} where Nazarov and Sodin gave a short exposition in terms of the standard metric in $\R^n$ instead of the \emph{canonical metric} $d(x,y)=(\E\brb{(F(x)-F(y))^2})^{1/2}$. The canonical metric has its own advantages in more general settings and it is clear how to compare $d(x,y)$ and $|x-y|$ but we think that is more natural to use the Euclidean metric for fields defined on domains in $\R^n$. 

The other way of thinking, that we call `analytic' emphasises how the field depends on $x$. Roughly speaking, we think that $F(x)$ is a function drawn at random from some space of functions. Alternatively, we can think that the probability measure is identified with a Gaussian measure on some function space. Probably the simplest way to construct a Gaussian function using this approach is to consider a random series
\[
F(x)=\sum a_i \phi_i(x)
\]
where $a_i$ are i.i.d. standard normal variables and $\{\phi_i(x)\}$ is a collection of linearly independent functions. Unless the collection is finite it is not immediately obvious that such series converges in any meaningful sense. Clearly, the answer depends on the collection of functions $\phi_i$. We will discuss this question a bit later.

On a more formal level, one can think that a Gaussian field is a map $\Omega\times X \rightarrow \R$, where $\Omega$ is the corresponding probability space. The `probabilistic' point of view emphasises the fact that for each $x\in X$ the maps $F(\cdot, x):\Omega\rightarrow \R$ is a Gaussian random variable. The `analytic' viewpoint emphasises properties of the function $F(\omega,\cdot):X\rightarrow \R$. 

These two descriptions are related and it should be possible to switch from one point of view to another depending on a particular question. Assuming that the series $F=\sum a_i \phi_i$ converges and we can swap summation and expectation, it is easy to see that 
\begin{equation}
\label{eq: K is sum}
K(x,y)=\E[F(x)F(y)]=\sum \phi_i(x)\phi_i(y).
\end{equation}

The other direction is less straightforward but still very well understood.  In this paper, we are interested only in the case $X=\R^n$ which significantly simplifies things. We also restrict our attention to the stationary case. This means that we want fields $F(\cdot)$ and $F(\cdot +z)$ to have the same distribution for every $z\in \R^n$. In terms of the covariance kernel $K$ it means 
\[
K(x,y)=K(x+z,y+z), \quad \forall x,y,z\in \R^n.
\]
In particular, this means that we can write $K(x,y)=K(x-y)$ where, abusing notations, we use $K$ for both functions of one and two variables. Since $K$ is a covariance kernel, it is a non-negative definite function. By Bochner's theorem  \cite[Theorem 5.4.1]{AdTa} there is a symmetric Borel measure $\rho$ (i.e. $\rho(A)=\rho(-A)$ for every measurable $A$) such that
\[
K(x)=\int_{\R^n} e^{ -2\pi i x\cdot t}\rho(d t).
\]

The measure $\rho$ is called the \emph{spectral measure} of the field $F$. Given this measure, we can construct a Hilbert space 
\[
L^2_{sym}(\rho)=\{\psi(t)\in L^2(\rho): \ \psi(-t)=\bar{\psi}(t)\}.
\]
This symmetry condition ensures that the corresponding Fourier transform
\[
\F(\psi)(x)=\int_{\R^n} \psi(t) e^{ -2\pi i x\cdot t}\rho(d t)
\]
is a real-valued function for every $\psi \in L^2_{sym}$. We define the Hilbert space 
\[
H=\F L^2_{sym}(\rho)=\left\{\int_{\R^n}e^{ -2\pi i x\cdot t}\psi(t)\rho(d t), \psi\in L^2_{sym}(\rho)\right\}
\]
where the scalar product is inherited from $L^2$. This is a reproducing kernel Hilbert space (RKHS) with the reproducing kernel $K$. Indeed, if $\phi=\F \psi$ then
\[
\begin{aligned}
\langle K_y(x),\phi(x)\rangle_H =&
\langle K(y-x),\phi(x)\rangle_H =
\langle e^{2\pi i y\cdot t}, \psi(t)\rangle_{L^2(\rho)}
\\
=&
\int_{\R^n} e^{-2\pi i y\cdot t} \psi(t) \rho(dt)=\phi(y).
\end{aligned}
\]

Given the Hilbert space $H$ one can construct the field as an isonormal process or white noise in $H$. We refer to \cite[Section 2.1]{NoPe} for more information about isonormal processes. Let $\{\psi_n\}$ be any orthonormal basis in $L^2_{sym}(\rho)$. We can formally define the white noise as 
\[
W=W_\rho=W_{L^2_{sym}(\rho)}=\sum_{k=1}^\infty a_k \psi_k(t),
\]
where $a_k$ are independent standard normal variables. This series obviously diverges in $L^2$, so one should treat it as a formal expression. If one wants to be rigorous, then the white noise is a collection of centred Gaussian random variables $Wh$ indexed by $h\in L^2_{sym}(\rho)$ such that
\[
\Cov(Wh_1,Wh_2)=\langle h_1, h_2\rangle_{L^2(\rho)}
=
\int_{\R^n} h_1(t)\bar{h}_2(t)\rho(dt).
\]
On the level of formal series, if $h_j=\sum_k b_{j,k} \psi_k$ for $j=1,2$, then 
\[
Wh_j=\langle W, h_j\rangle=\sum a_k b_{j,k}, \qquad j=1,2
\]
and 
\[
\begin{aligned}
\Cov(Wh_1,Wh_2)=&\E \left[\sum a_k b_{1,k}\sum a_k b_{2,k}\right]
\\
=&\sum b_{1,k}b_{2,k}=\langle h_1,h_2\rangle.
\end{aligned}
\]
Using this language we can formally write
\[
F(x)=(\F W)(x)=W e^{-2\pi i x\cdot t }=\int_{\R^n} e^{-2\pi i x\cdot t } \sum a_k \psi_k(t) \rho(dt)= \sum a_k \phi_k(x),
\] 
where $\phi_k=\F\psi_k$ form an orthonormal basis in $H$. As in \eqref{eq: K is sum} we can see that $F$ is a Gaussian field with the covariance kernel $K(x,y)$. In particular, we can see that it does not depend on our choice of the orthonormal basis. 

There is a particularly important case when the spectral measure $\rho$ is continuous with respect to the Lebesgue measure. In this case, we can write $\rho(dt)=\rho(t)dt$ and it is easy to see that
\[
W_\rho=\sqrt{\rho(t)}W_{dt}.
\]
In particular, in this case, we can rewrite everything in terms of the standard $L^2$ space and obtain the following representation:
\[
F(x)=W*q=Wq(x-\cdot)
\]
where $*$ is the convolution,  $W$ is the canonical white noise in $L^2$ and 
\[
q=\int \exp(-2\pi i x \cdot t) \sqrt{\rho(t)}dt
\]
is a convolution square root of $K$ i.e. $K=q*q$. This representation is also known as the \emph{moving averages} representation \cite[Section 5.3]{AdTa}.

Let us conclude this section by emphasising one more time that explanations in this section are correct only on the `physics level of rigour', but they could be written absolutely rigorously.

\subsection{Percolation}

In this section, we give a \emph{very} brief introduction to percolation theory. A proper introduction would require writing a whole book.  Our goal is to introduce the main notations and terms that we will be using throughout this survey. More information about percolation can be found in the classical book of Grimmett \cite{Grimmett}, Bollob{\'a}s and Riordan present a different take on the subject that also includes some results about the scaling limit of the critical percolation. There are also good lecture notes by Werner \cite{Werner_lecture_notes} and Duminil-Copin \cite{DC_lecture_notes}. 

There are two main versions of the Bernoulli percolation: \emph{bond} percolation and \emph{site} percolation. In both cases, we start with a graph. For the purpose of this survey we restrict the attention to regular lattices in $\R^2$. Two main examples are the square lattice and the triangular lattice. We choose some parameter $p$ and independently declare edges (bond percolation) or vertices (site percolation) open with probability $p$ or closed with probability $1-p$. Percolation clusters are connected subgraphs made of open bonds or sites. In the context of fields, the most relevant is the bond percolation. 

\begin{figure}
\includegraphics[width=0.45\textwidth]{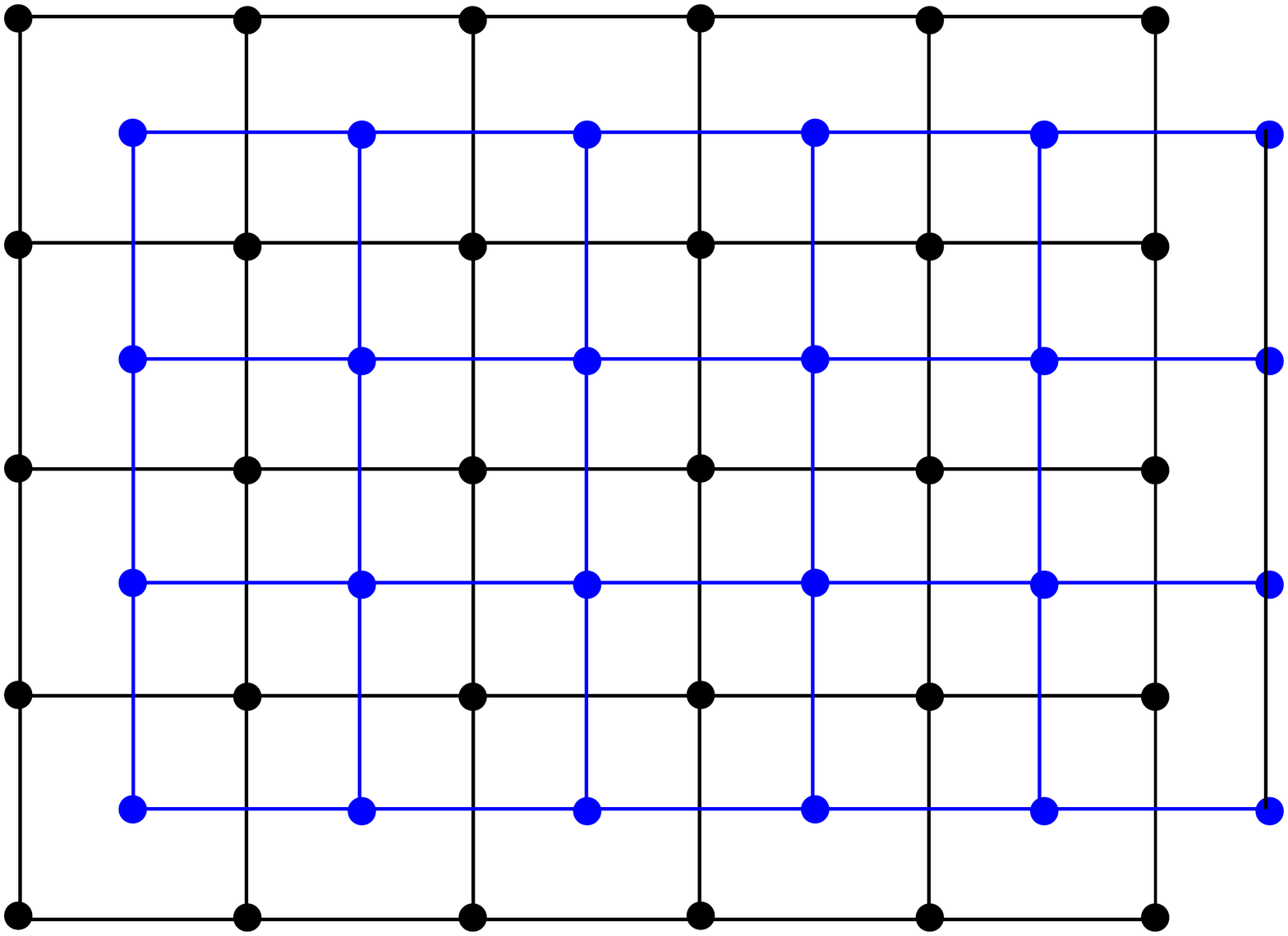}
\hspace{0.05\textwidth}
\includegraphics[width=0.45\textwidth]{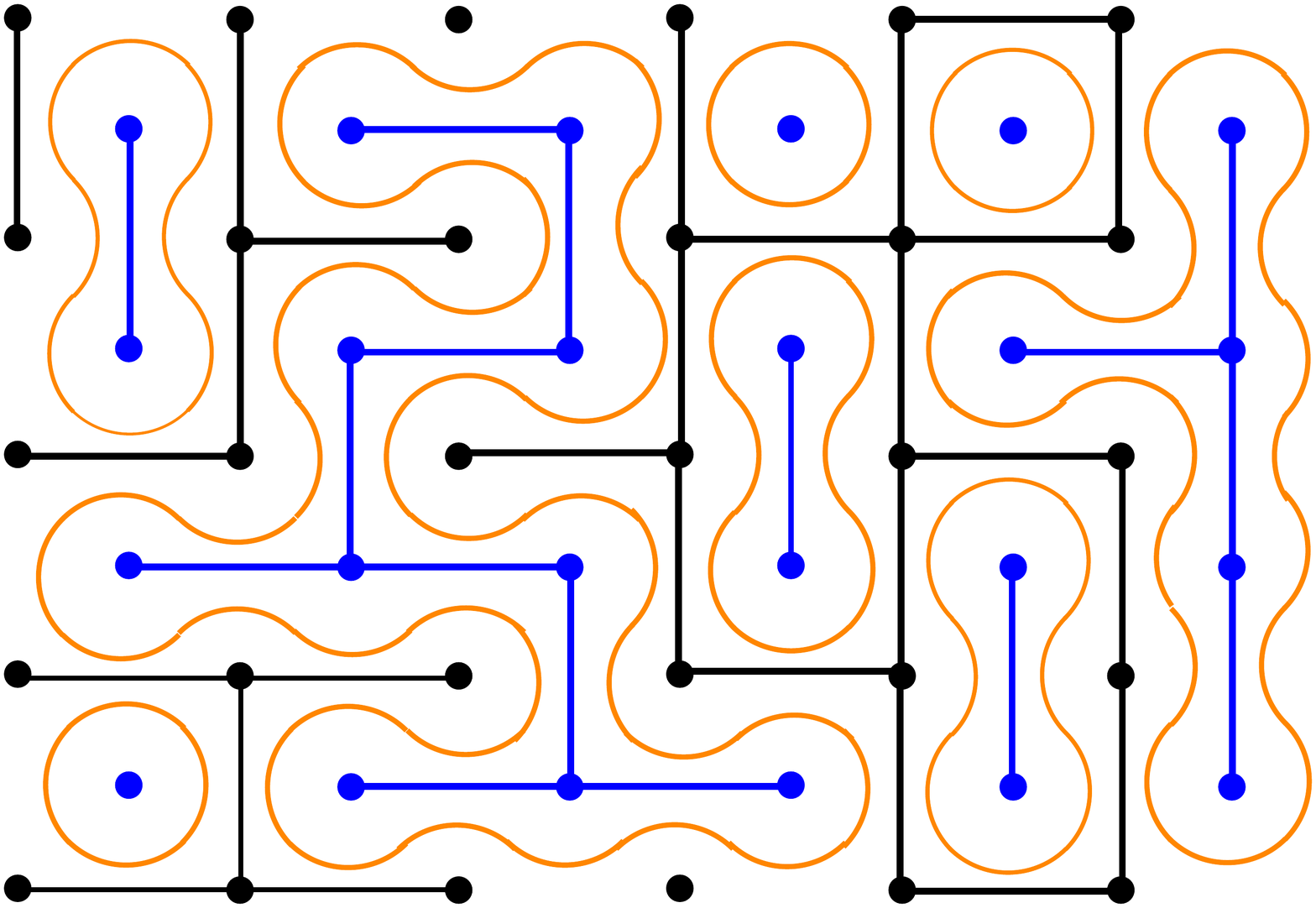}
\caption{Left: a part of the square lattice and its dual. Right: an example of the critical percolation ($p=1/2$) on both lattices and curves separating primary and dual clusters.}
\label{fig: percolation}
\end{figure}

Quite often it is also convenient to consider simultaneously the graph and its dual. When we declare a bond in the primary graph to be open, we declare the dual bond closed and vice versa. This dual picture is particularly nice for the square lattice since its dual graph is also the square lattice. See Figure \ref{fig: percolation} for an example. 

One of the first questions about percolation models is whether there is an infinite cluster. We define $\theta(p)$ to be the probability that a fixed vertex, say the origin, belongs to an infinite cluster. The critical probability is
\[
p_c=\inf_{p\in[0,1]}\{p: \theta(p)>0\}.
\]
A priori, it is not obvious that for a given graph the critical $p_c$ is non-trivial i.e. $0<p_c<1$. It is also not immediately clear whether there is an infinite cluster at $p=p_c$ i.e whether $\theta(p_c)=0$ or not. 

For the square lattice $\Z^2$ the critical probability $p_c=1/2$. The simplest way to argue that this must be the case is by self-duality. One can argue that at $p=1/2$ percolations on the primal and dual lattices are the same, hence at $1/2$ the model is self-dual and so it must be critical. There are several other  heuristics giving $p_c=1/2$. 
\begin{theorem}[Kesten \cite{Kesten}]
\label{thm: Kesten}
For the Bernoulli percolation on $\Z^2$ the critical probability $p_c=1/2$. For $p\le 1/2$ there is no infinite cluster with probability one, for $p>1/2$ there is an infinite cluster with probability one. 
\end{theorem}
By now there are several rigorous proofs of this result. The first proof is due to Kesten. A simpler and shorter proof was given by Bollob{\'a}s and Riordan in \cite{BoRi_short_Kesten}. 

The critical probability is lattice dependent, but it is believed that under mild assumptions on the regularity of the lattice the large scale behaviour at criticality is universal i.e. lattice independent. It is also conjectured that critical percolation has a conformally invariant scaling limit. In 1999 Schramm \cite{Schramm} introduced a one-parameter family of random curves defined in terms of the Loewner Evolution. They are called Stochastic or Schramm Loewner Evolution (SLE) curves. It is conjectured that interfaces between clusters (see Figure \ref{fig: percolation}) converge to SLE(6) curves. Since then, the proof has been simplified  \cite{Beffara_cardy,KhSm} and the entire argument is very short and elegant. It has been extended to the convergence of the complete collection of the percolation interfaces \cite{CaNe} to a generalization of SLE that is called the \emph{Conformal Loop Ensemble} or CLE. So far all proofs rely heavily on the structure of the triangular lattice and all attempts to generalize it to different lattices were futile. Some partial results were obtained by Beffara\cite{Beffara_mesoscopic}. See also a discussion of universality in \cite{Beffara_universality}.

\section{Main examples}

There are two main examples that we should keep in mind: the Random Plane Wave (RPW) and the Bargmann-Fock field. Here we discuss their two-dimensional versions, but they have natural higher-dimensional generalizations. 

\subsection{Random Plane Wave}

Let us start with the formal definition. The RPW field can be defined as a stationary Gaussian field in $\R^2$ with the covariance kernel $K(x)=J_0(|x|)$, where $J_0$ is the $0$-th Bessel function. From this description we can see that point-wise covariances go to zero as points move away from each other, but the rate is quite slow since 
\[
J_0(t)=\sqrt{\frac{2}{\pi t}}\cos\left(t-\frac{\pi}{4}\right)+O(t^{-3/2}). 
\]
Another important observation is that $K$ is an oscillating function.

The corresponding spectral measure is the normalized arc-length on the unit circle. This is a very singular measure, which is yet another indication that covariance is slowly decaying. On the other hand, it is compactly supported, hence the field is real analytic. Since the support is on the unit circle, the field is an eigenfunction of the Laplacian. If one thinks naively of the white noise in $L^2_{sym}(\rho)$ as a collection of independent Gaussians at each point of the unit circle, then one can think that the field is a random linear combination of plane waves with the same energy travelling in all possible directions. In this sense, one can think that the RPW is a natural notion of a `typical' eigenfunction of the Laplacian in $\R^2$. 

\begin{figure}
\includegraphics[width=0.45\textwidth]{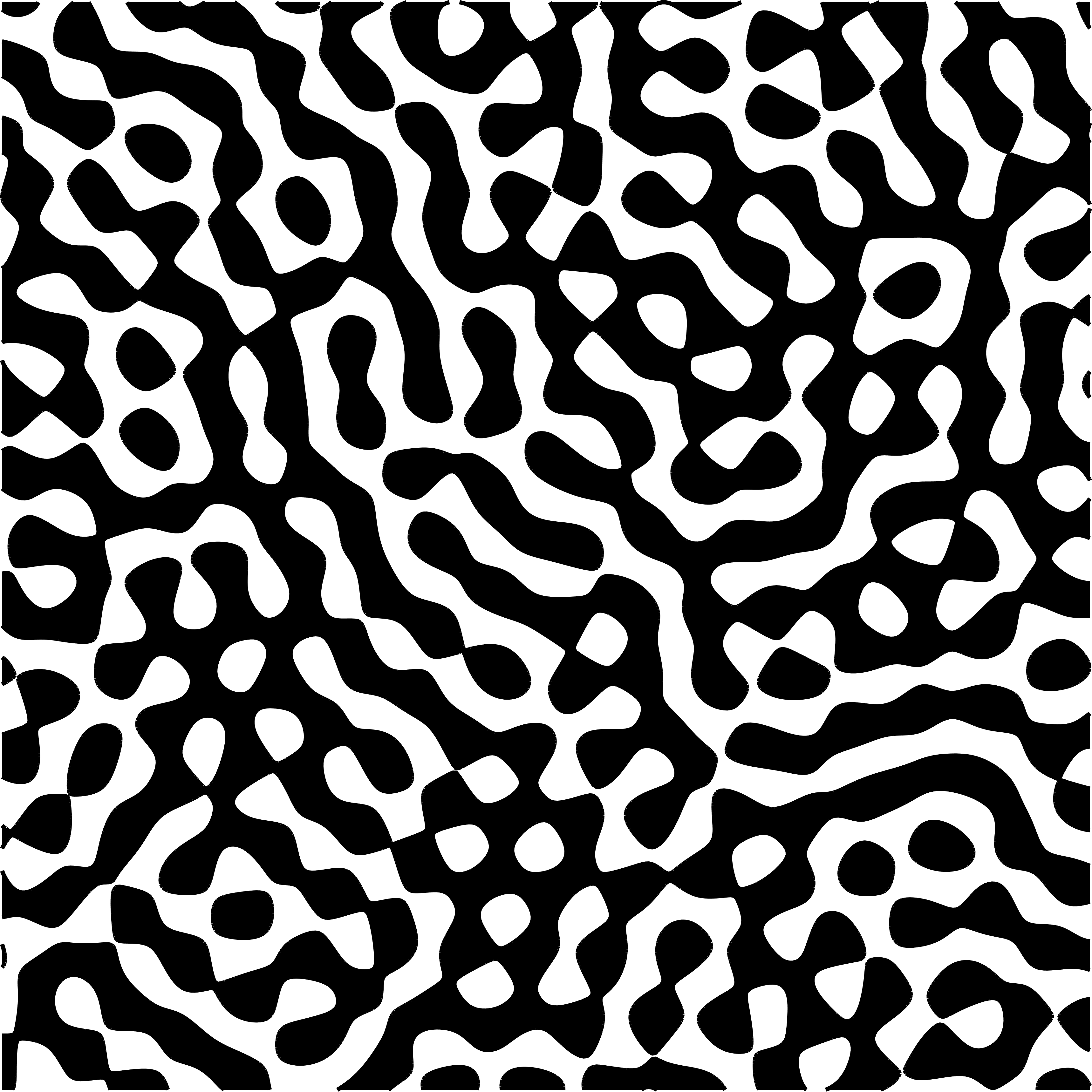}
\hspace{0.05\textwidth}
\includegraphics[width=0.45\textwidth]{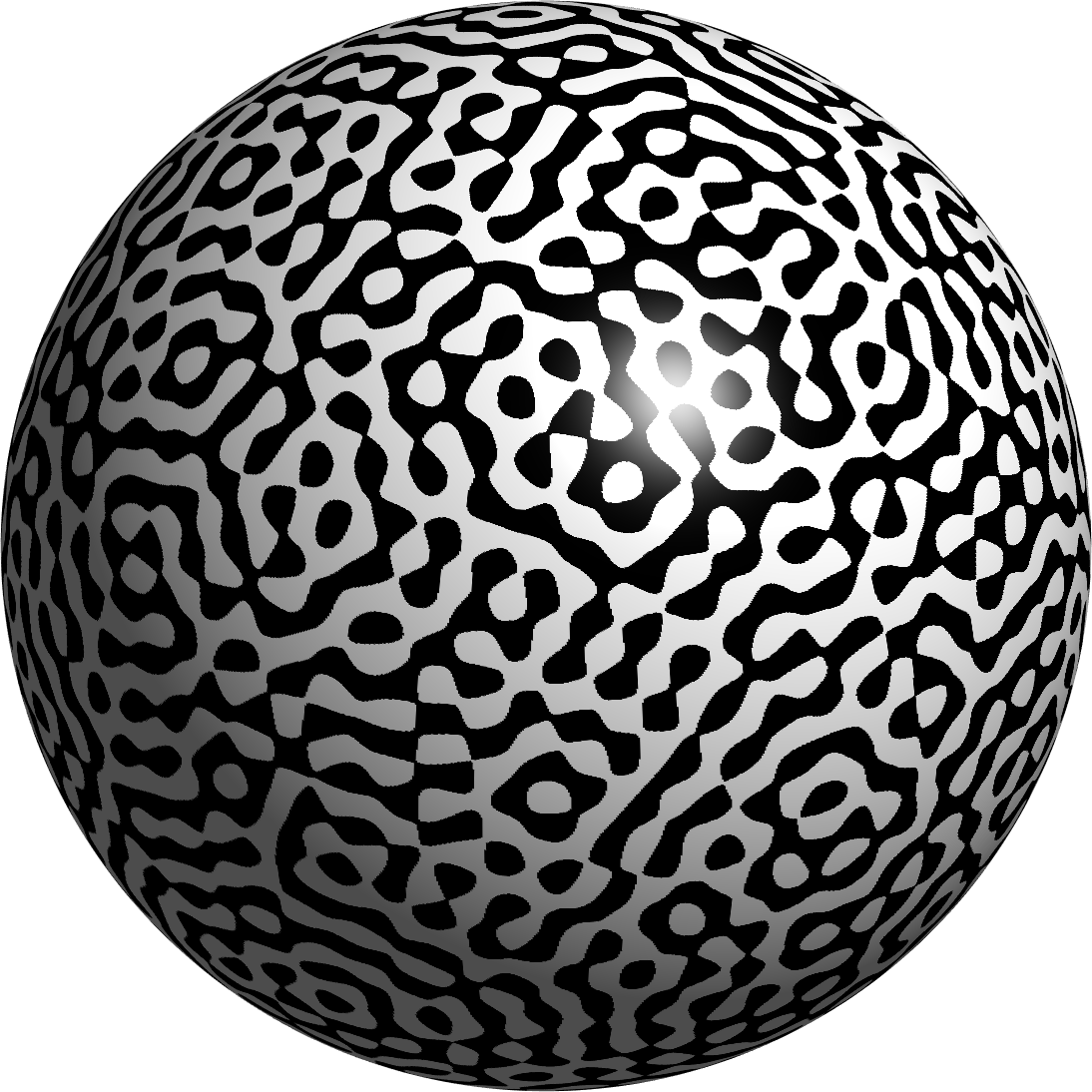}
\caption{The random plane wave (left) and the random spherical harmonic of degree $100$ (right). Black and white are regions where the fields are negative or positive.}
\label{fig: RPW and RSH}
\end{figure}

Finally, the RPW has the following series representation  
\begin{equation}
\label{eq: RPW series}
F(x)=\sum_{-\infty}^\infty C_n J_{|n|}(|x|)e^{in\theta},
\end{equation}
where $C_k$ are independent \emph{complex} Gaussian random variables subject to the condition $C_{-k}=\bar{C}_k$. See Figure \ref{fig: RPW and RSH} (left) for a sample of the random plane wave.

This field also appears as the scaling limit of `narrow-band' functions \cite{Zelditch}. In the case of the sphere, eigenfunctions of the Laplacian are spherical harmonics. For a fixed degree $n$ the space of spherical harmonics is of dimension $2n+1$ and the corresponding eigenvalue is $-n(n+1)$. The eigenspace is equipped with the natural $L^2$ metric. The random spherical harmonic of degree $n$ is defined as the standard Gaussian vector in the eigenspace. To be more precise, if we choose any orthonormal basis, for example, the standard $Y_n^k$, then the random spherical harmonic of degree $n$ is 
\[
F_n(\theta, \phi)=\sum_{k=-n}^{n}a_n^k Y_n^k(\theta,\phi),
\] 
where $(\theta,\phi)$ are the standard spherical coordinates and $a_n^k$ are independent $N(0,1/(2n+1))$ random variables. With this normalization the expected $L^2$ norm of $F$ is $1$. Figure \ref{fig: RPW and RSH} (right) shows an example of a random spherical harmonic.

From the addition formula for the Legendre polynomials one can see that the covariance kernel of this field is
\[
K_n(x,y)=\E[ F_n(x)F_n(y)]=\frac{1}{2n+1}\sum_{k=-n}^{n}Y_n^k(x)Y_n^k(y)=P_n(\cos \theta(x,y)),
\]
where $x,y$ are two points on the unit sphere, $P_n$ is the Legendre polynomial of degree $n$ and $\theta(x,y)$ is the angle (i.e spherical distance) between $x$ and $y$. 

From this formula one can immediately see that this is a (spherically) stationary field. Alternatively, this could be seen from the invariance of the eigenspace and $L^2$ norm under rotations. 

The random spherical harmonics converge to the random plane wave as $n\to \infty$. To be more precise, let $x_0$ be any point on the sphere and $\exp_{x_0}:T_{x_0}S^2\rightarrow S^2$ be the exponential map from the tangent plane to the sphere. Then we can define
\[
f_n(x)=F_n(\exp_{x_0}(x/n)).
\]
Note that we rescale by $n$ which is approximately the square root of the eigenvalue. This means that we rescale so that the wavelength becomes of order $1$. Clearly, the covariance of this field is
\[
P_n(\cos \theta(\exp_{x_0}(x/n),\exp{x_0}(y/n))).
\]
Since the exponential map is almost isometry near $x_0$, the spherical distance between images is almost the distance between the points, hence uniformly in $x$ and $y$ the covariance behaves as
\[
P_n\left(\cos \frac{|x-y|}{n}\right)\to J_0(|x-y|),
\]
where the limit follows from Hilb's asymptotics for Legendre polynomials.

It is possible to define a similar field on any compact manifold. In the generic case there are no large eigenspaces like in the spherical case. Instead, we take $L^2$ normalised eigenfunctions of the Laplacian $\phi_k$ and consider 
\[
f_n(x)=\sum_{k=n^2}^{n^2+n} a_k \phi_k(x)
\]
where $a_k$ are i.i.d $N(0,1/n)$. That is, like in the spherical case, we take a linear combination of order $n$ eigenfunctions around $n^2$-th eigenfunction. It is possible to show, see \cite{Zelditch} that under some mild assumptions of the manifold, a similarly rescaled field converges to the random plane wave. 

\subsection{Bargmann-Fock Feild}

\begin{figure}
\includegraphics[width=0.45\textwidth]{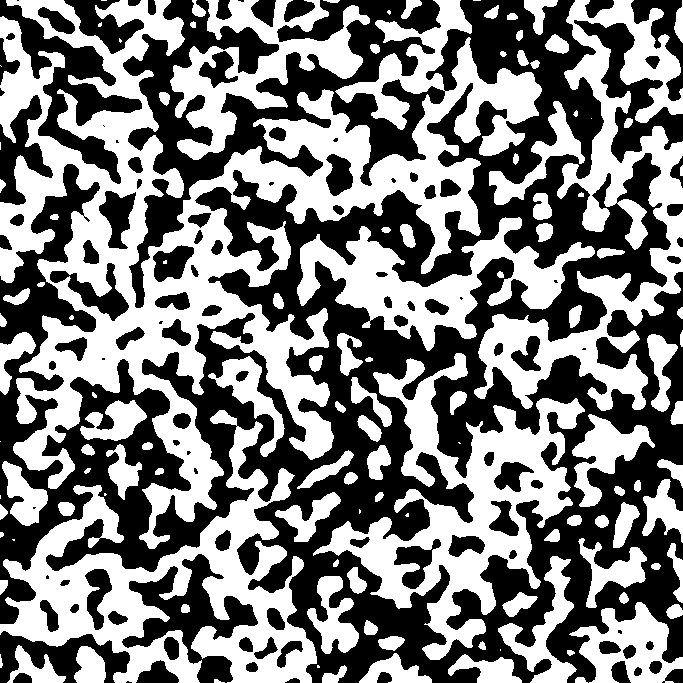}
\hspace{0.05\textwidth}
\includegraphics[width=0.45\textwidth]{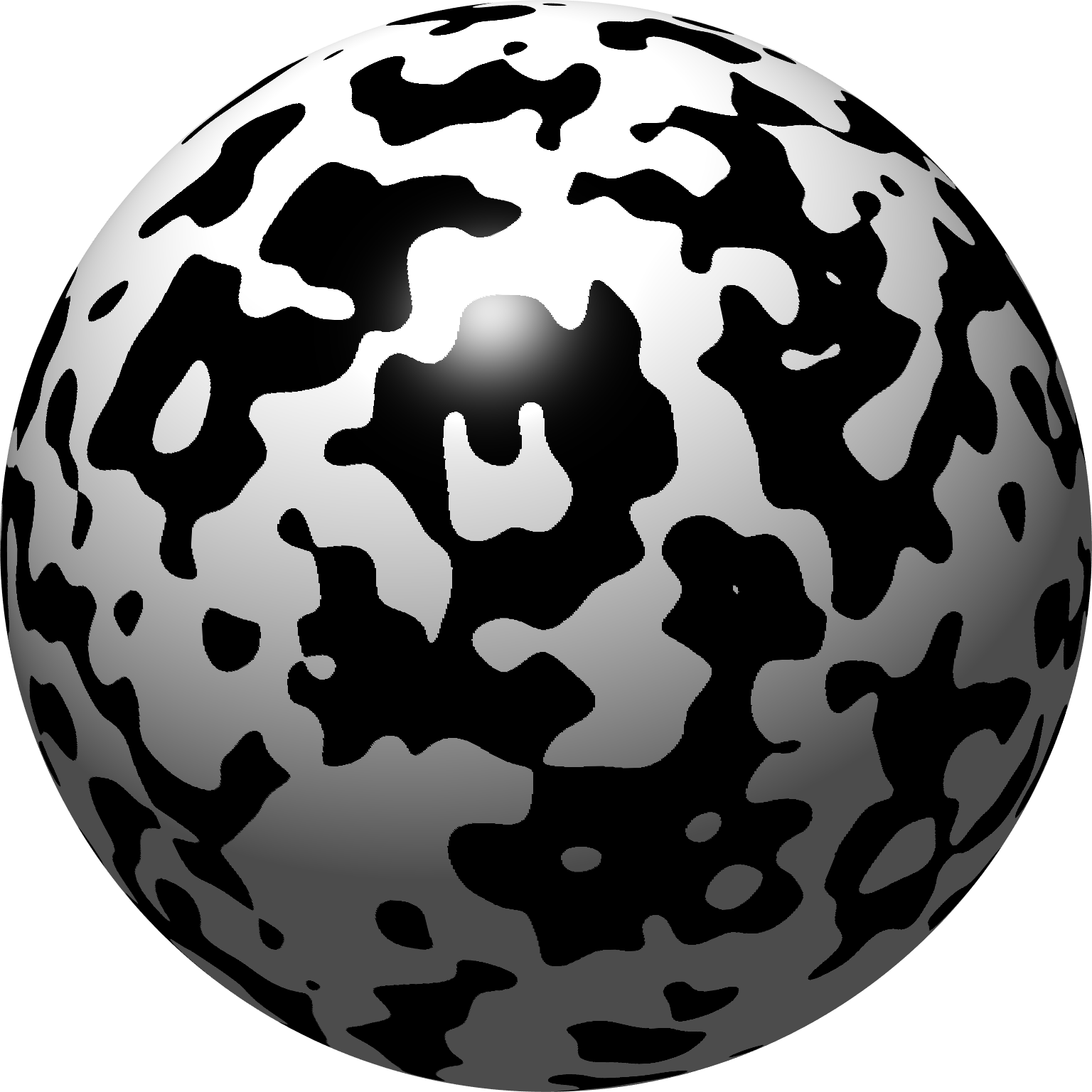}
\caption{Bargmann-Fock (left) and Kostlan ensemble of degree $300$ (right). Black and white are regions where the fields are negative or positive.}
\label{fig: BF and Kostlan}
\end{figure}

We define the Bargmann-Fock field in $\R^2$ as
\begin{equation}
\label{eq:BF def}
F(x)=F(x_1,x_2)=e^{-|x|^2/2}\sum_{m,n=0}^\infty a_{n,m}\frac{1}{\sqrt{n!m!}}x_1^nx_2^m.
\end{equation}
It is easy to see that the covariance kernel is
\[
\begin{aligned}
\E\brb{ F(x)F(y)}=&e^{-|x|^2/2}e^{-|y|^2/2}\sum \frac{(x_1y_1)^n(x_2y_2)^m}{n!m!}
\\
=&e^{-|x|^2/2}e^{-|y|^2/2}e^{x\cdot y}=e^{-|x-y|^2/2}.
\end{aligned}
\]
Since the covariance depends on $|x-y|$ only, the field is stationary and isotropic. Since the kernel is the Gaussian density function, the spectral measure has also Gaussian-type density. It is possible to give a more direct description of the corresponding Hilbert space, but at least formally, it could be described as the only Hilbert space for which functions in \eqref{eq:BF def} form an orthonormal basis. This, space is very closely related to the space of holomorphic functions known as Bargmann, Segal-Bargmann or Bargmann-Fock space that appears in quantum mechanics \cite{Bargmann}.  See Figure \ref{fig: BF and Kostlan} for a sample of the Bargmann-Fock field.

Similarly to the RPW, the Bargmann-Fock field is a scaling limit of a certain family of fields on the unit sphere (see the right plot in Figure \ref{fig: BF and Kostlan}). Let $H_n$ be the space of all homogeneous polynomials of degree $n$ in $\R^3$. We equip it with the scalar product such that the following functions form an orthonormal basis
\[
\sqrt{\frac{n!}{\alpha_1!\alpha_2!\alpha_3!}}x^\alpha,
\]
where $x^\alpha=x_1^{\alpha_1}x_2^{\alpha_2}x_3^{\alpha_3}$ and $|\alpha|=\alpha_1+\alpha_2+\alpha_3=n$.
Up to a constant factor the corresponding scalar product is the same as the Bargmann-Fock product 
\[
\langle f_1, f_2\rangle=\int_{\C^3}f_1(z)\bar{f}_2(z)e^{-\|z\|^2}dz.
\]
This is the real trace of the only structure on complex polynomials that is invariant under the unitary group (although there are many real ensembles that are invariant with respect to the orthogonal transformations \cite{Kostlan}). Thus it is possible to claim that this is the most natural Hilbert space structure on $H_n$. 

The scalar product on $H_n$ gives rise to a canonical Gaussian measure. The corresponding field
\[
F_n(x)=\sum_{|\alpha|=n}a_\alpha \sqrt{\frac{n!}{\alpha_1!\alpha_2!\alpha_3!}}x^\alpha
\]
where $a_\alpha$ are independent standard Gaussian random variables. This field is known as the Kostlan ensemble or real Fubini-Study ensemble. Direct computation shows that the covariance kernel for this field is
\[
K_n(x,y)=\cos^n \theta(x,y).
\]

As before, we can use the exponential map to project this field to a tangent plane
\[
f_n(x)=F_n(\exp_{x_0}(x/\sqrt{n})).
\]
Similarly, we can compute the asymptotic behaviour of the covariance
\[
\cos^n \theta(\exp_{x_0}(x/\sqrt{n}),\exp_{x_0}(y/\sqrt{n}))\approx \cos^n \left(\frac{|x-y|}{\sqrt{n}}\right)\to e^{-|x-y|^2/2}.
\]
This implies that the scaling limit of the Kostlan ensemble is the Bargmann-Fock field.

\begin{remark}
As we have explained above, the Kostlan ensemble is the natural candidate for a model of a `typical' homogeneous polynomial of degree $n$. Hence, its zero or nodal lines are the natural model of a `typical' real projective variety of degree $n$. Thus the results about these nodal lines could be interpreted as the results about `typical' projective varieties.
 \label{rem: Kostlan projective variety}
\end{remark}
 
\section{Number of level sets}

We can think that a Gaussian field describes a landscape where $f(x)$ is the height function. Using this analogy, it is very natural to represent the field via the collection of its level sets $\{f=\ell\}$ or excursion sets $\{f\ge \ell\}$. This is exactly how we use isolines on topographic maps. Under mild regularity assumptions a typical level set is a collection of smooth hypersurfaces since $f$ and $\nabla f$ can not vanish simultaneously. This is known as Bulinskaya's lemma. It was originally proved by Bulinskaya in the one-dimensional case \cite{Bulinskaya}, see \cite[Lemma 11.2.10]{AdTa} for a proof of the general statement.

\subsection{One-dimensional case}
The one-dimensional case is particularly simple. In this case, level sets are just isolated points and excursion sets are intervals between them. In the context of Gaussian processes, the level sets are known as \emph{level crossings}. Quite often they are separated into up- and down-crossings. The study of level sets, especially nodal sets, where $\ell=0$ has a long history and goes back to the works of Kac \cite{Kac} and Rice \cite{Rice}. One of the main advantages is that the number of level crossings is a \emph{local} observable. This means that the number of level crossings in a disjoint union of two sets is equal to the sum of level crossings inside these sets. In particular, this allows writing an integral formula for moments of the number of level crossings. This formula is now known as the Kac-Rice formula. It holds in much higher generality, but we formulate only the simplest version. A detailed discussion of this formula and its proof in a very general setting can be found in \cite[Sections 11.2 and 11.5]{AdTa}.

\begin{theorem}[Kac-Rice Formula] 
Let $f$ be a Gaussian process on an interval $I$. We assume that $f$ has $C^1$-paths, that it $f\in C^1$ with probability one. We also assume that for pairwise distinct points $t_1,\dots,t_k\in I$ the joint distribution of $(f(t_1),\dots,f(t_k))$ is non-degenerate. Let $N_\ell$ be the number of points where $f(t)=\ell$, then the $k$th  factorial moment 
\[
\E[N_\ell^{[k]}]=E[N_\ell(N_\ell-1)\cdots(N_\ell-k+1)]
\]
is given by the following multiple integral
\begin{equation}
\label{eq:Kac-Rice}
\iint_{I^k}\E\brb{\prod_{i=1}^k|f'(t_i)|\mid f(t_1)=\cdots=f(t_k)=\ell}p_{t_1,\dots,t_k}(\ell,\dots,\ell)\prod_{i=1}^k dt_i,
\end{equation}
where  $p_{t_1,\dots,t_k}$ is the joint density of $(f(t_1),\dots,f(t_k))$. 

In particular, 
\[
\E\brb{N_\ell}=\int_I\E\brb{|f'(t)|\mid f(t)=\ell}p_{t}(\ell)dt.
\]
\end{theorem}

One of the nice things about this formula is that for Gaussian processes the conditional expectation can be computed in terms of the covariance kernel. Indeed,
\[
(f'(t_1),\dots,f'(t_k),f(t_1),\dots,f(t_k))
\]
is a Gaussian vector. Its covariance can be written explicitly in terms of the covariance kernel of $f$ and its derivatives. Using the Gaussian regression formula we can explicitly compute the conditional expectation. Although in theory it can be done, in practice, even with $k=2$ the resulting expression could be quite complicated and the analysis of the multiple integral becomes non-trivial. 

There is one case when this formula becomes particularly simple. Let us assume that $f$ is a centred stationary process with covariance kernel $K$, then the law of $(f'(t),f(t))$ does not depend on $t$. It is centred and its covariance matrix is
\[
\begin{pmatrix}
-K''(0)& 0\\ 0& K(0)
\end{pmatrix}.
\]
We can immediately see that $f'(0)$ and $f(0)$ are independent.
So the integral becomes 
\[
\begin{aligned}
|I|\E\brb{|f'(0)|}p_{f(0)}(\ell)=&|I|\sqrt{\frac{-2K''(0)}{\pi}}\frac{1}{\sqrt{2\pi K(0)}}e^{-\ell^2/2}
\\
=&\frac{|I|}{\pi}\sqrt{\frac{-K''(0)}{K(0)}}e^{-\ell^2/2}.
\end{aligned}
\]

\subsection{Higher dimensions}
In higher dimensions the situation is much more complicated. First of all, the level sets now can have non-trivial geometry and topology. Moreover, even a simple observable like the number of the connected components of a level set is more complicated. The main difference is that it is no longer an additive quantity, so there is no integral formula similar to the Kac-Rice formula to compute it. One can still write a Kac-Rice formula in a higher dimensional case, but it will not count the number of components, instead, it will compute the volume of the level set. Since this is not the main focus of our survey, we will not discuss it at length, but a lot is known about additive quantities like volume. Their expected volume can be easily computed in the stationary case. For some fields it is even possible to write good estimates on the variance and even prove (central) limit theorems. Since it is only tangentially related to our main question about the large scale behaviour of level sets, we will not discuss this line of research any further. Instead, we refer to a survey by Wigman \cite{Wigman_survey} and lecture notes by Berzin, Latour and Le\'on \cite{BLL} and references wherein. 

For a long time the question about the expected number of level sets was completely out of reach. The breakthrough came in 2009 when Nazarov and Sodin made an observation that the number of level sets is semi-local, that is the number of `large' level sets is small. Initially, they proved the law of large numbers for random spherical harmonics. 
\begin{theorem}[Nazarov-Sodin \cite{NaSo_RSH}]
\label{thm: NS RSH}
Let $f$ be the random spherical harmonic of degree $n$. Then there exists $c>0$ such that for every $\epsilon>0$ we have 
\[
\P\brb{\left|\frac{N(f)}{n^2}-c \right|>\epsilon}\le C(\epsilon)e^{-c(\epsilon)n},
\]
where $N(f)$ is the number of nodal domains and $c(\epsilon)$ and $C(\epsilon)$ are constants that depend on $\epsilon$ only.
\end{theorem}
In particular, this strong concentration theorem implies the law of large numbers, namely that $N(f)/n^2\to a$ a.s.

Later, they extended this result to stationary fields in $\R^n$ that satisfy rather mild assumptions:
\begin{enumerate}
\item The spectral measure has finite $4$-th moment. This assumption implies that the field almost surely is $C^{2-\epsilon}$ for every $\epsilon>0$.
\item The spectral measure has no atoms. This implies ergodicity with respect to the shifts.
\item The spectral measure is not supported on a linear hyperplane. This implies that the distribution of $\nabla f$ is non-degenerate.
\end{enumerate}

\begin{theorem}[Nazarov-Sodin \cite{NaSo_asymptotic}]
\label{thm: NS stationary}
Let $f$ be a stationary field in $\R^n$ such that its spectral measure satisfies the three assumptions above. Let $S$ be a bounded open convex set containing the origin and $N_S(R,f)$ be the number of nodal domains of $f$ that are contained in the rescaled set $R\cdot S$, then there is $c \ge 0$ such that
\[
\frac{N_S(R,f)}{R^n \mathrm{Vol}(S)}\to c.
\]
The convergence is a.s. and in $L^1$. Moreover, under additional mild assumptions $c>0$. 
\end{theorem}

\begin{remark}
Although the theorem above was stated for the number of nodal domains, the same proof works for the number of level/excursion sets for all levels. 
\end{remark}

It is also possible to estimate the variance of the number of level/excursion sets. There are two complementary results. In \cite{NaSo_variance} Nazarov and Sodin linked nodal lines of Gaussian fields to a percolation model on a certain degree four graph and obtained a low bound on the number of nodal domains from the variance of the number of percolation clusters. Their argument could be viewed as a step towards justification of the Bogomolny-Schmit heuristics. As a result, they have shown that if $f$ is a $C^3$ stationary field and its covariance kernel decays polynomially with some positive exponent, then the variance of the number of nodal domains in a ball of radius $R$ grows faster than $R^\sigma$ for some positive $\sigma$. On one hand, the conclusion is rather weak since $\sigma$ is not specified, on the other hand, the assumptions of the field are absolutely minimal. In particular, the result is applicable to the RPW. 

At about the same time we obtained a complementary result. In \cite{BMM_variance} we have shown that if $f$ satisfies some regularity assumptions, its covariance is non-negative and decays sufficiently fast and the spectral density near the origin is bounded away near the origin, then for a generic level the variance of the number of level sets inside a ball of radius $R$ grows at least as $R^2$. We refer to \cite[Theorem 2.7]{BMM_variance} for the precise list of assumptions and the precise explanation for which levels it holds. We would like to point out that the proof \emph{does not hold} for $\ell=0$, our assumptions are much stronger than that of \cite{NaSo_variance} but our conclusion is also much stronger and it is sharp for a wide class of fields. 

Imposing even stronger assumption it is  possible to prove a Central Limit Theorem that is valid for all levels and in all dimensions \cite{BMM_CLT}. 
\begin{theorem} Let $f$ be a stationary field in $\R^d$  with the covariance $K=q*q$. Assume that $q\in C^5$, $q$ and its derivatives decay faster that $|x|^\beta$ for some $\beta>9d$, then the variance of the number of level/excursion sets $N(R,\ell)$ inside the cube $\Lambda_R=[-R,R]^d$ grows like the volume
\[
\frac{\mathrm{Var} \brb{N(R,\ell)}}{\mathrm{Vol}(\Lambda_R)}\to \sigma^2, \qquad R\to\infty
\]
and
\[
\frac{N(R,\ell)-\E\brb{N(R,\ell)}}{\sqrt{\mathrm{Vol}(\Lambda_R)}}\to Z ,
\]
where $Z$ is a normal random variable with variance $\sigma$. Moreover, if $\int K>0$, then $\sigma>0$ for all $\ell$. 
\end{theorem}

We do not believe that our assumptions are necessary. For example, the natural bound on all moments of the number of critical points would allow to reduce the decay rate to $\beta>3d$. There is no reason to believe that it is the best possible result, most probably, the result should hold for $\beta>0$. Similarly, we do not think that $C^5$ assumption is optimal.
 
\section{Percolation Conjecture}

\subsection{Number of nodal domains for the random plane wave}
The laws of large numbers discussed above led naturally to the question about the nature of the constant $c$ mentioned above. Sometimes this constant is referred to as the \emph{Nazarov-Sodin} constant, but it is important to remember that it is really a functional that depends on the spectral measure and the level. 

The question about the number of nodal domains also appears in quantum chaos. Physicists  Blum,  Gnutzmann and Smilansky \cite{BGS} proposed that the nodal domain statistics could be used as a criterion for quantum chaos. Same year Bogomolny and Schmit \cite{BSch02} proposed a percolation model for the number of nodal domains of the random plane wave. Later in \cite{BSch07} they used (non-rigorous) Harris criterion \cite{Harris_criterion} to strengthen the percolation conjecture. In \cite{BDSch07} Bogomolny, Dubertrand and Schmit gave further numerical support to the percolation conjecture	 by showing numerically that nodal domains of the random plane waves converge to SLE(6).

Roughly speaking in \cite{BSch02} Bogomolny and Schmit have argued that on average nodal lines form a square lattice, but since nodal lines do not intersect almost surely, these intersections should be resolved one way or another. By symmetry, both possibilities are equally likely, so the nodal structure could be identified with bond percolation on the square lattice with $p=1/2$. In the case of the square lattice, the number of percolation clusters per vertex is known explicitly. This gives a prediction for the number of nodal domains per unit volume. This is essentially the constant appearing in Theorems \ref{thm: NS RSH} and \ref{thm: NS stationary}. To be more precise, they have conjectured that the expected number of nodal domains in a rescaled domain $R\cdot S$  is 
\[
\frac{\mathrm{Area}(S)R^2}{2\pi}\frac{3\sqrt{3}-5}{\pi}\approx \frac{\mathrm{Area}(S)R^2}{2\pi} 0.0624373.
\]
Moreover, they also claimed that the variance of the number of nodal domains also scales as $R^2$. Their conjecture was supported by numerical simulations. 

Later, more careful computations \cite{Nastasescu11, Konrad12, BeKe} have shown that the number $0.06243$ is not quite right, the real value is approximately $0.0589$. This shows that the prediction of Bogomolny and Schmit is not quite right. On the other hand, the number of percolation clusters per vertex is not a universal quantity, it is lattice-dependent. Nodal lines are clearly not a small perturbation of the square lattice, so there is no reason to believe that the square lattice percolation would have exactly the same density of clusters. Moreover, there is one more scaling in this argument. To compare the number of clusters per vertex and the number of domains per unit volume, we have to find the right relation between vertices and volume. Bogomolny-Schmit proposed to relate the wavelength and the lattice mesh. One can also argue that it is equally reasonable to relate the number of vertices and the number of critical points (after all lattice faces play the role of local minima and maxima and bond intersections are related to saddle points). This scaling gives a close, but a different prediction. We would go even further and claim that one of the really interesting questions is why the Bogomolny-Schmit prediction is off just by $5\%$. 

It is possible to interpret nodal domains of the random plane wave in terms of a percolation-type model. 
\begin{figure}
\includegraphics[width=0.31\textwidth]{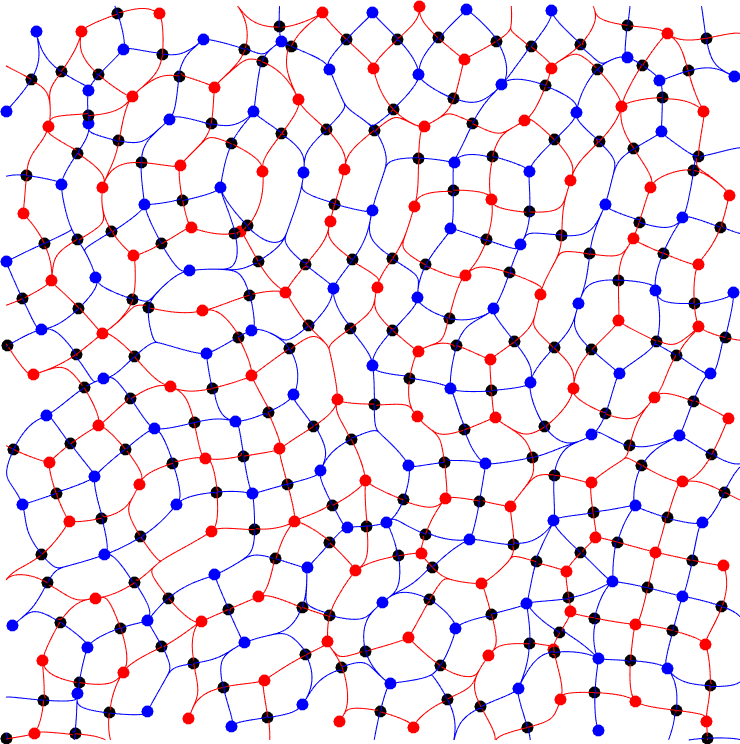}\hspace{0.02\textwidth}
\includegraphics[width=0.31\textwidth]{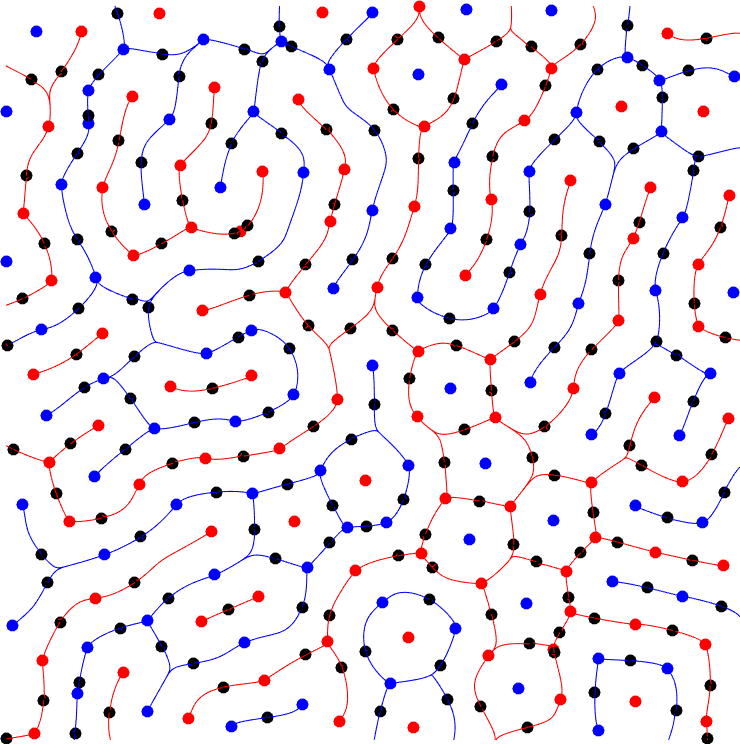}\hspace{0.02\textwidth}
\includegraphics[width=0.31\textwidth]{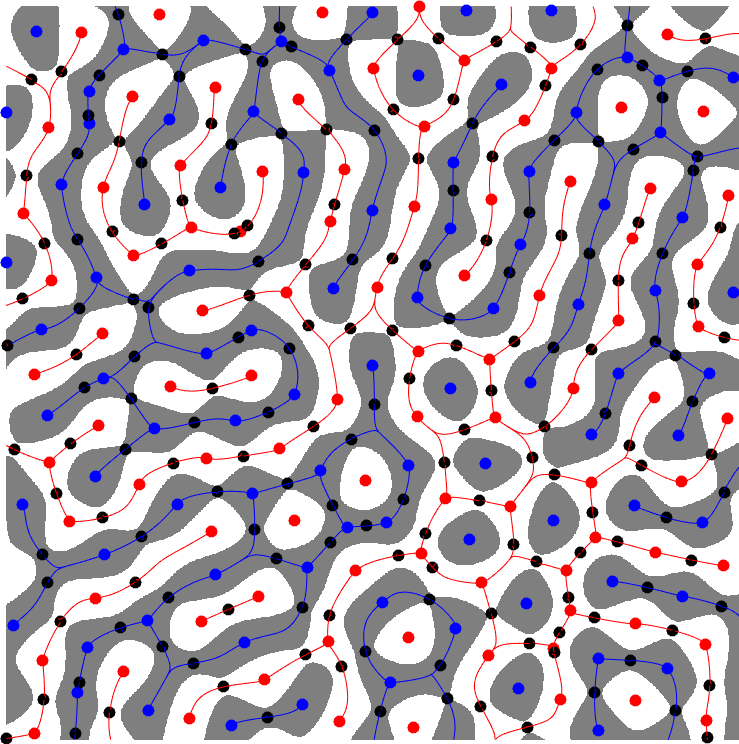}
\caption{Graph of flow-lines (left), corresponding percolation (center) and nodal domains (right).}
\label{fig: flow lines}
\end{figure}
Since the random plane wave is an eigenfunction of the Laplacian, all local minima are negative and local maxima are positive. Let us consider gradient flow lines  from a saddle that go along eigenvectors of the Hessian. Two ascending lines will terminate at two local maxima and two descending lines will terminate at local minima. This way we have two planar graphs. The `positive' graph has its vertices at local maxima and two maxima are connected by a bond if there is a saddle such that its ascending flow lines terminate at these maxima. Similarly, we define the `negative' graph. These two graphs are dual to each other and intersections of their bonds are saddles. See Figure \ref{fig: flow lines} for a sample of these graphs. We do the following percolation on these graphs: if the field is positive at a saddle, then we keep the `positive' (or ascending) bond and remove the negative one, if it is negative, then we keep the `negative' bond and remove the `positive' one. Clearly, by symmetry, we choose one or the other with the probability $1/2$. It is easy to see that the resulting percolation cluster corresponds exactly to nodal domains. 

On one hand, this is a bond percolation-type model, but it is a horrible one. First of all, we do percolation on a random graph, secondly, signs of the functions at different saddles are not independent (moreover, correlation decays slowly since the Bessel function decays slowly) and finally, how we choose bonds is correlated with the graph itself. 

So far there was no progress in analysing this percolation model. On the other hand, it is plausible, that this picture explains why Bogomolny-Schmit prediction is so close to the true value. Since the `positive' and `negative' graphs have the same distribution, they are self-dual. In particular, on average a lot of their parameters are the same as that of the square lattice. For example, the average face is a quadrangle and the average degree is $4$. Moreover, simulations suggest that the probability that a vertex has a large degree has very light tails. Most of the percolation clusters are small, so the number of percolation clusters is mostly defined by the local structure of the graph. Probably this is the reason why cluster density for the square lattice and for the flow-line graph are so close, but it is not clear whether it is possible to make rigorous sense out of this hand-waving argument. 
 
\subsection{Universal observable}
As we have explained above, the number of percolation clusters per vertex is lattice-dependent and hence it is very hard to get it right. On the other hand, $p=1/2$ is the critical probability for the bond percolation on the square lattice. It is conjectured that all `global' observable in critical percolation models are universal, i.e. independent of the lattice. The percolation interpretation of the nodal domains is clearly self-dual, so it is very natural to assume that it is critical as well. It is quite natural to re-interpret the claim of Bogomolny and Schmit and to conjecture that the nodal domains of the random plane wave are in the same `universality class' as the critical percolation. In particular, to conjecture that many global observables have the same behaviour. 

In the critical percolation context, one of the most natural observables is the crossing probability. That is the probability that in an $n\times \lambda n$ rectangle there is an open cluster connecting the left side to the right side. It is conjectured that this probability has a scaling limit as $n\to \infty$ and the limit does not depend on the lattice. The limiting probability is given by an explicit formula in terms of $\lambda$. This formula was first predicted by Cardy \cite{Cardy92} using the Conformal Field Theory. It has been proven only in the case of the critical percolation on the triangular lattice \cite{Smirnov01}.

The same observable can be defined for the RPW or any other smooth field in the plane. We ask whether there is a positive nodal domain connecting the left and right sides of an $n\times \lambda n$ rectangle. Equivalently, this is a probability that there is a curve crossing the rectangle such that the field is positive along this curve. This crossing probability for the random plane wave has been numerically computed in \cite{BeKe} and it matches Cardy's formula very well. This supports the following form of the Bogomolny-Schmit conjecture

\begin{conjecture}[Bogomolny-Schmit conjecture]
All large scale connectivity properties of the RPW nodal lines and domains are the same as for the critical percolation. In particular, all crossing events have the same scaling limits. The collection of all nodal lines has a scaling limit which is conformally invariant and given by the Conformal Loop Ensemble CLE($6$). 
\end{conjecture}  

It is also very natural to conjecture that the relationship between the RPW and percolation goes beyond the nodal domains. If we consider some other level $\ell\ne 0$, then the large scale properties of excursion sets should correspond to the percolation with non-critical $p$. We call this statement the \emph{off-critical Bogomolny-Schmit conjecture}.

These two conjectures are still wide open and there was almost no progress in understanding the connection between the RPW and percolation. 
There are some numerical evidence and heuristics in support of this conjecture but there are no rigorous results supporting it. On the other hand, there is one result that is somewhat contrary to the off-critical conjecture. In \cite{BMM_variance} it was shown that for non-zero $\ell$ the variance of the number of excursion sets in a square of size $R$ grows at least as $R^3$ whereas for percolation the variance of the number of clusters grows as $R^2$ both in critical and off-critical case. It is important to note that the number of clusters is non-universal (although the variance growth rate is) and is not directly related to the large scale connectivity properties, so this result does not immediately disprove the percolation conjecture, but it shows that \emph{some} large scale observables behaves differently. 

\subsection{Generalized Bogomolny-Schmit conjecture}

Over the last decade it became clear that the connection between the RPW and percolation is not unique to the RPW, it seems to be true for a wide class of smooth Gaussian fields. It is still not clear what is exactly the class of fields for which the correspondence should hold, but we believe that some statement along the following lines should hold. 

\begin{conjecture}
\label{conj: generalized BS}
Let $f$ be a centred stationary Gaussian field in $\R^2$ satisfying the following assumptions
\begin{itemize}
\item With probability one $f$ is sufficiently smooth. 
\item The distribution of $f$ and its derivatives is non-degenerate.
\item The covariance kernel is sufficiently symmetric.
\item The covariance kernel decays sufficiently fast.
\item The covariance kernel is (strongly) positive.
\end{itemize}
Then the nodal domains of $f$ should have the same large scale properties as critical percolation clusters and excursion sets for non-zero levels should behave like off-critical percolation clusters.
\end{conjecture} 

Some smoothness is clearly needed at least to guarantee that level sets are smooth curves. A priori, it is not clear why the joint distribution of $f$ and its derivatives at various points should be non-degenerate for the conjecture to hold, but so far all partial results require some form of non-degeneracy. On the other hand, it is not a very strong restriction. For example, if the support of the spectral measure contains an open set, then one can usually prove the required non-degeneracy without too much effort. It is more complicated with fields like the RPW that have singular spectral measure.

Some symmetry is also clearly necessary, without it there are obvious counterexamples. The exact requirement is not obvious. Currently, most of results are valid under the assumption that $K(x_1,x_2)=K(x_1,-x_2)$ which is probably very close to the minimal requirement.

The decay of correlation is very natural: percolation is completely local, different bonds are open or closed completely independent, so for the field it is natural to require that distant parts are weakly correlated. What is the minimal decay of correlation is not clear, but but there are reasons to believe that kernel decaying faster than $|x|^{-\beta}$ for $\beta>2$ is enough. This requirement will be discussed later.

So far, many results require either that $K\ge 0$ of that $K=q*q$ with $q\ge 0$ (strong positivity). The last assumption is also far from obvious and it is not obvious whether it is required. In particular, the covariance kernel of the RPW is not positive.  This assumption is needed in many partial results since it guaranteed that the field is positively associated (see Theorem \ref{thm Pitt} below) which allows using many percolation techniques. The strong positivity assumption works very well with the moving averages representation $f=W*q$. 

The fact that nodal lines should correspond to the critical percolation is not really surprising. By symmetry of the normal distribution, $f$ and $-f$ have the same distribution. Hence the distributions of $\{f\le 0\}$ and $\{f\ge 0\}$ are the same. This means that the model is self-dual, so it should be critical. This is very similar to the heuristics implying that $p_c=1/2$ for the bond percolation on $\Z^2$ (Kesten's theorem \ref{thm: Kesten}). 

Finally, we would like to point out that from this point of view the Bogomolny-Schmit conjecture for the random plane wave looks even more surprising. As we have mentioned above, its kernel is not positive, it decays rather slowly and some of its derivatives are deterministically related since it is an eigenfunction of the Laplacian. On the other hand, the Bargmann-Fock field is real analytic, the joint distributions of various values and derivatives are non-degenerate, and the covariance kernel is positive and decays super-exponentially fast. 

Let us conclude this section with two important non-examples. 

\begin{example}
Let us consider the field that is a version of the discrete white noise. We tile the plane by unit squares. The field is constant within each square and values at different squares are independent standard Gaussians. This field is not stationary, but if we randomly shift tiling it will be stationary. The nodal domains of this field are exactly the same as site percolation clusters on the square lattice with $p=1/2$ which is \emph{not critical}. So in this case nodal domains do not behave like critical percolation.
\end{example}

\begin{example}
The other example is arguable the most well-studied Gaussian field: the Gaussian Free Field. This is a $\log$-correlated field. Since the covariance kernel blows up on the diagonal, the field is not even a function, it is a random distribution. Despite it being a distribution, there is a way to make sense of its nodal lines and it has been shown \cite{SchSh} that these nodal lines are Schramm-Loewner Evolution SLE($4$) curves while the percolation interfaces are conjectured to converge to SLE($6$) curves, so in this case the nodal lines \emph{are conformally invariant} but belong to a different universality class.   
\end{example}

\subsection{Previous results related to percolation}

Long before the work of Bogomolny and Schmit there were several results about percolating properties of the level sets of Gaussian fields. To be more precise, one can ask whether the level/excursion sets of a Gaussian field are almost surely bounded or not. 

Probably the first result in this direction is a series of works of Molchanov and Stepanov \cite{MSI,MSII,MSIII}. They were looking at a wide class of fields both on $\Z^d$ and $\R^d$, not necessary Gaussian, and discussed when the set $\{f\le \ell\}$ contains almost surely a unique unbounded component. In particular, in \cite[Theorem 3.5]{MSII} they have shown that if $f$ is a continuous stationary Gaussian field with the covariance kernel decaying sufficiently fast (essentially faster than $|x|^{-d}$), then there is a finite level $\ell$ such that $\{f\le \ell\}$ percolates a.s. 

The second important result that we want to mention is a more recent one by Alexander. Theorem 3.4 of \cite{Alexander} states that if $f$ is a stationary $C^1$ Gaussian field in $\R^2$ such that its covariance kernel is non-negative and decays to zero at infinity, then all level lines are almost surely bounded. In other words level lines never percolate. 

Finally, we would like to mention one result that is seemingly unrelated.  We say that random variables $x_i$ are positively associated if 
\[
\Cov\brb{f(x_1,\dots,x_n),g(x_1,\dots,x_n)}\ge 0
\]
for all functions $f$ and $g$  that are non-decreasing in each variable. In \cite{Pitt} Pitt have proved the following theorem: 
\begin{theorem}
\label{thm Pitt}
Jointly Gaussian random variables $x_i$ are positively associated if and only if $\Cov\brb{x_i,x_j}\ge 0$ for all $i$ and $j$. 
\end{theorem}
In other words, positively correlated Gaussian variables are positively associated. Despite its simplicity, this is a surprisingly non-trivial statement. 

In the context of Gaussian fields we have the notion of increasing events. These are the events such that if it holds for a field $f$, then it holds for all larger fields. In particular, all connectivity events for excursion sets are increasing events. Pitt's theorem implies that if the covariance kernel is non-negative, then all increasing events are positively correlated. Pitt's theorem is stated for finite-dimensional Gaussian vectors, but it is not very hard to show that standard approximation techniques extend it to the field setting (see \cite[Lemma A12]{RiVa} for details).

 This result is a field equivalent of the Harris \cite{Harris} or Fortuin, Kasteleyn, Ginibre (FKG) \cite{FKG} inequalities that give positive correlations of increasing events in percolation. Such inequalities play a very important role in percolation theory since they allow to combine various crossing events in a controlled way. Although there are some approaches to the percolation that do not rely on Harris-FKG, it is still a fundamental part of the theory. This is the main reason why the positive kernel is one of the assumptions in our generalised Bogomolny-Schmit conjecture \ref{conj: generalized BS}.

\section{Many small steps towards the percolation conjecture}

\subsection{Locality}

One of the defining features of percolation models is locality: different bonds are open or closed independently. In particular, this means that percolation pictures within disjoint domains are independent. Locality plays an important role in the percolation theory and in the study of its scaling limit. SLE($6$) was initially identified as the (conjectured) scaling limit of percolation interfaces because it is the only SLE that has some form of the locality property. 

Naively, the locality in random fields should follow from the decay of correlation, but the question is more subtle. Let's consider two disjoint domains that are far apart. Although point-wise values are weakly correlated, there are uncountably many points, so collectively, they could be strongly correlated. 

Indeed, if we look at our main examples, they are real analytic functions, so if we know the field in any open set, then we know it everywhere. Equivalently, the sigma-algebra generated by the restriction of the field to any open set is the full sigma-algebra. This means that level lines in disjoint domains can not be really independent. Fortunately, it turns out that relevant events, such as crossing events, indeed decorrelate when domains are well separated. 

The relevant notion in the case of fields is quasi-independence of crossing and other similar events. For crossing events we would like to obtain a statement of this kind: let $B_1$ and $B_2$ be two disjoint rectangles, then the covariance between crossing events in $R\cdot B_1$ and $R\cdot B_2$ goes to $0$ as $R$ goes to infinity. Sometimes we also need to control the rate of convergence. More generally, we would like to have a similar estimate for all events describing the topology of nodal level sets. 

Below we will describe a range of such results. In this section we always implicitly assume that all fields are centred, stationary, sufficiently smooth (usually $C^3$ is enough) and non-degenerate (different theorems use different notions of non-degeneracy, but in all examples it is usually easy to check these conditions). Below we will only concentrate on the more important assumption on the decay of the covariance kernel. 

One of the first results in this direction was obtained by Nazarov,  Sodin and Volberg in \cite[Theorem 3.2]{NSV07} (see also \cite[Theorem 3.1]{NS11}). They worked in the context of the Gaussian entire function (GEF) which is complex-valued rather than real-valued, but its correlations also decay fast similarly to the Bargmann-Fock field. They have shown that the restriction of GEF to two well separated compacts could be written as restrictions of two independent copies plus a small perturbation. Essentially the same argument would work for the Bargmann-Fock and other fields with very fast decaying covariance. This is not exactly what we need, one also has to show that the small perturbation does not change the structure of the nodal set with high probability.  

The first, very natural approach to such results is discretization. The main idea is rather straightforward. Let us evaluate the field inside the box $R\cdot B$ on a lattice with mesh $\epsilon$. This gives us a \emph{finite} Gaussian vector with dimension of order $R^2\epsilon^{-2}$. If the mesh is fine enough, then it should be possible to argue that knowledge of the field in these points captures all relevant information with very high probability. For example, if we know the field on a mesh then we know whether nodal domains percolate or not with very high probability. Now if we have two boxes, then we have two finite-dimensional vectors. If boxes are far apart, then correlations between each component of one vector and each component of the other vector are small. 

One of the first results in this direction is due to Mischaikow and Wanner \cite{MiWa} who used discretization of the field to compute the homology of the nodal set. A better result was obtained by Beffara and Gayet in \cite{BeGa}. They were studying the crossing events (we will discuss it in more detail a bit later) and used exactly the approach mentioned above. They have shown \cite[Theorem 1.6]{BeGa} that if $\epsilon=o(R^{-8-\delta})$ for some small $\delta>0$, then with probability almost one the nodal structure of the discretized field is the same as the structure of the original field. After that, they gave an estimate of the dependence of Gaussian vectors. 
\begin{theorem}[Theorem 4.3 \cite{BeGa}]
Let $X$ and $Y$ be two Gaussian vectors of dimensions $n$ and $m$ correspondingly. We consider two joint distributions $\mu$ and $\nu$ such that under $\mu$ they are independent and under $\nu$ pairwise covariances $\Cov(X_i,Y_j)\le \delta$ for all $i$ and $j$, then we have the following bound on the total variation distance:
\begin{equation}
\label{eq: BG total variation}
d_{\mathrm{TV}}(\mu, \nu)\le \mathrm{const}(m+n)^{8/5}\delta^{1/5}.
\end{equation}

\end{theorem} 
This means that two correlated vectors are essentially independent as long as $\delta=o((m+n)^{-8})$. Combining this with an estimate on the mesh size we get that nodal structures in two domains of size $R$ that are $R$ apart are essentially independent if the covariance kernel decays faster than $|x|^{-144}$. 

Both of these estimates can be improved. In \cite[Theorem 1]{BeMu_discr} we have shown that discretization gives a faithful representation as long as $\epsilon=o(R^{-2-\delta})$ and in \cite[Proposition 16]{BeMu_discr} we have improved an estimate \eqref{eq: BG total variation} by showing that the covariance of all events depending only on signs of $X$ and $Y$ correspondingly (all events describing nodal sets are of this type)  is bounded by
\begin{equation}
\label{eq: BM total variation}
 \mathrm{const}(m+n)^{4/3}\delta^{1/3}\log^{1/3}(1/\delta).
\end{equation}

The best to day result in this direction was obtained by Rivera and Vanneuville \cite[Theorem 1.12]{RiVa}. There they have shown that crossing events are quasi-independent if the covariance kernel $K(x)=o(|x|^{-4})$. Their argument is quite different. Although they do use discretization as an intermediate step, their estimates are \emph{uniform} in the lattice step.

Inspired by the work of Rivera and Vanneuville as well as earlier work of Piterbarg \cite{Piterbarg} we have obtained in \cite{BMR} an explicit formula for the covariance of a certain class of `topological' events. Unlike all previous results, our approach is dimension-independent so we are able to obtain quasi-independence in higher dimensions. We also do not assume that the field is centred, hence, without loss of generality, we consider nodal sets instead of general level sets. 

Before explaining the main result and the main ideas behind the proof, we have to introduce some notations. To simplify the explanation, we do it in the simplest possible setting and we will skip some technicalities. Full details can be found in the original publication. 

Let $B=[a_1,b_1]\times\dots\times[a_n,b_n]$ be a rectangular box in $\R^n$. We treat it as a stratified $n$-dimensional manifold. We decompose it as the union of its $n$-dimensional interior ($n$-dimensional strata), open $(n-1)$-dimensional sides etc up to all corners that are $0$-dimensional strata. We will be considering isotopies that preserve this stratification, that is isotopies of $B$ that send each stratum to itself. We say that two compact subsets of $B$ are equivalent if there is a stratified isotopy of $B$ that sends one set to another. Elementary topological events are events that the nodal sets/domains of a field inside $B$ belong to a certain equivalence class. Topological events are events from the $\sigma$-algebra generated by elementary topological events.  This is a rigorous way of saying that the nodal set has a certain structure. This notion on one hand allows enough flexibility so that the analytic rigidity is not a problem anymore. On the other hand, this $\sigma$-algebra is rich enough to describe all events that we are interested in such as the number of components or various crossing events. 

\begin{figure}
\includegraphics[width=0.8\textwidth]{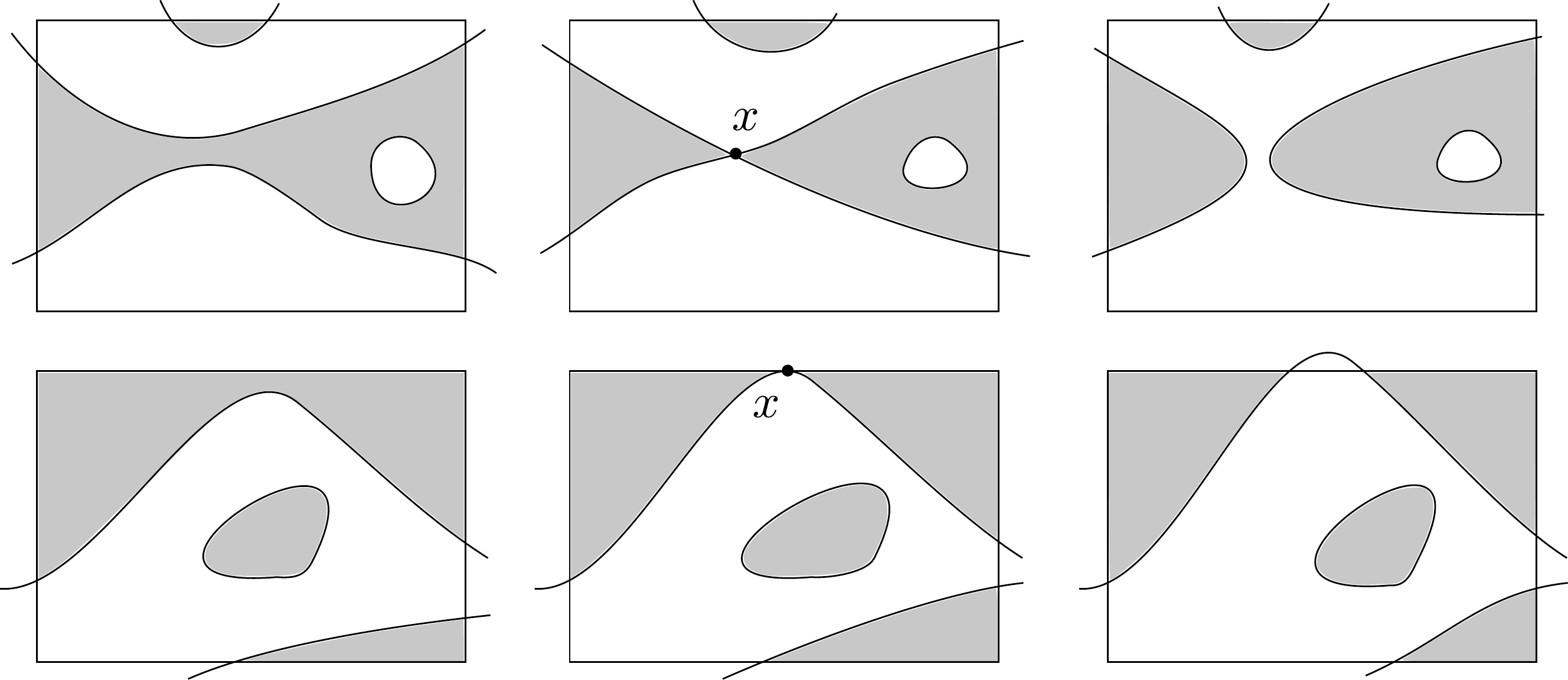}
\caption{Central panels show functions with pivotal points for the lef-right crossing event. Left and right show small perturbations. In the top row the pivotal point is a saddle inside the domain, in the bottom row the pivotal point is an extremum of the field restricted to a one-dimensional stratum.}
\label{fig:pivotal}
\end{figure}

Next, we need a field equivalent of the notion of a \emph{pivotal} point in percolation. Let us fix some topological event $A$. We say that a point $x$ is \emph{positively pivotal} for a function $f$ if small positive local perturbations of $f$ near $x$ belong to $A$ and small negative local perturbations do not and we say that it is \emph{negatively pivotal} if it is the other way round. One should think that $f$ is on the boundary of a topological event and positivity/negativity give us a notion of `orientation'. Note that not all boundary functions have pivotal points. Up to a small exceptional set of boundary functions pivotal points are non-degenerate stratified critical points of $f$ with $f(x)=0$. See Figure \ref{fig:pivotal} for an illustration. 

Let $f$ be a Gaussian field with the covariance kernel $K(x,y)$\footnote{We do not assume that $f$ is stationary}. We assume that $f\in C^2$ and non-degenerate. Let $B_1$ and $B_2$ be two boxes and $A_1$, $A_2$ be two topological events in these boxes. Then there are two positive \emph{pivotal measures} $\pi^\pm$ such that
\begin{equation}
\label{eq: topological covariance}
\P[A_1\cap A_2]-\P[A_1]\P[A_2]=\int_{B_1\times B_2} K(x,y)(d\pi^+(x,y)-d\pi^-(x,y)).
\end{equation}

These measures are also stratified. For each pair of strata in $B_1$ and $B_2$ these measures have a component that is continuous with respect to the product of the corresponding Lebesgue measures on these strata. The density at $(x,y)$ is very closely related to the joint density of two stratified critical points at $x$ and $y$ subject that they are pivotal for $A_1$ and $A_2$ correspondingly. The positive measure $\pi^+$ corresponds to the case when either both points are positively pivotal or negatively pivotal and the negative measure $\pi^-$ corresponds to the case when one is positively pivotal and the other is negatively pivotal. These densities are essentially the same as the conditional expectations that appear in Kac-Rice formula \ref{eq:Kac-Rice}. 

The precise formula for the densities is rather long, complicated and involves an interpolation between a pair of identical copies of $f$ and a pair of independent copies. All details can be found in \cite{BMR}.

This formula has a couple of important corollaries. When both events are increasing, there are no negatively pivotal points, hence the negative term $\pi^-$ disappears. If we also assume that $K\ge 0$ then the integrand in \eqref{eq: topological covariance} is non-negative and we recover Pitt's positive association result. Densities appearing in the definition of $\pi^\pm$ could be uniformly bounded. If $B_i$ are two boxes of the size of order $R$ then the covariance can be uniformly bounded by a term comparable to
\[
R^{2d}\sup_{x\in B_1,y\in B_2}K(x,y).
\] 
This means that if the covariance kernel decays is $o(|x-y|^{-2d})$ then all topological events are quasi-independent. When $d=2$ this is exactly the quasi-independence result from \cite{RiVa}. 

Finally, \eqref{eq: topological covariance} leads to some form of Harris criterion or rather to a necessary condition. Let $R\cdot B_1$ and $R\cdot B_2$ be two disjoint boxes rescaled by a large factor $R$. If we assume the `percolation hypothesis', then, as $R$ tends to infinity, crossing events in the rescaled boxes should decorrelate. Pivotal densities appearing in the definition of $\pi^+$ (there is no $\pi^-$ since events are increasing) are similar to densities of `four-arm saddles at scale $R$' i.e. saddle points $x$ such that four level lines $\{f=f(x)\}$ stay separated within a ball of radius $R$. So the covariance of crossing events should be comparable to
\[
I(R)^2\int_{R\cdot B_1\times R\cdot B_2}K(x-y)dxdy,
\]
where $I(R)$ is the density of four-arm saddles at scale $R$. These `four-arm' saddles are the field equivalent of `four-arm' events in the percolation theory. Assuming the percolation hypothesis $I(R)$ should behave like $R^{-\zeta_4}$, where $\zeta_4$ is the so-called \emph{four-arm exponent}. By Kesten scaling relations \cite{Kesten_scaling} $\zeta_4=d-1/\nu$ where $\nu$ is the \emph{correlation length exponent}. Roughly speaking it is defined so that a Bernoulli percolation with non-critical $p$ is indistinguishable from the critical percolation on scale $R$ as long as $|p-p_c|\ll R^{-1/\nu}$. So, assuming the percolation hypothesis, we should have
\[
R^{2/\nu-2d}\int_{R\cdot B_1\times R\cdot B_2}K(x-y)dxdy\to 0.
\]
This is almost the same criterion as the field version of the Harris criterion that appeared in \cite{BSch07} (see also \cite{Weinrib}). 

It is believed that in dimension $d=2$ the correlation length exponent $\nu=4/3$ for all lattices. So far, it has been proven only for the triangular lattice \cite{SmWe}. Assuming this exponent, the Harris criterion becomes 
\begin{equation}
\label{eq: Harris criterion}
R^{-5/2}\int_{R\cdot B_1\times R\cdot B_2}K(x-y)dxdy\to 0.
\end{equation} 
Roughly speaking, for positive $K$ this corresponds to $K$ decaying faster than $|x|^{-3/2}$. Note that the covariance kernel for the random plane wave is decaying like $|x|^{-1/2}$, but it is also \emph{oscillating}, so it could be argued that it satisfies this form of the Harris criterion. This gives a (rather weak) support to the percolation conjecture for the random plane wave. 

Although this is not completely rigorous and is only a necessary condition, the fact that it is essentially the same as the Harris criterion derived by Bogomolny and Schmit using a completely different approach gives some support to this heuristic.

\subsection{Russo-Seymour-Welsh inequality}

Crossing probability plays a very important role in percolation theory. Proving Cardy's formula was the most important step in Smirnov's proof of the convergence of critical percolation interfaces on the triangular lattice to SLE($6$). So far, we are very far from proving Cardy's formula for any Gaussian field. 

The first step in understanding the scaling properties of percolation is the Russo-Seymour-Welsh  (RSW) inequality that was first obtained independently by Russo \cite{Russo} and by Seymour and Welsh \cite{SeWe}. 
\begin{theorem}[Russo-Seymour-Welsh]
Let us consider the Bernoulli percolation with $p=1/2$ on the square lattice. Then for every $\rho>1$ there is $c>0$ which depends on $\rho$ but not on $n$ such that
\[
c<\P\brb{\text{there is a horizontal crossing in } [0,\rho n]\times [0,n]}<1-c,
\]
\end{theorem}
Similar statements have been proved for other percolation models as well. We use the term RSW for all such statements. Other, shorter and more general proofs were obtained by Bollob{\'a}s and Riordan in \cite[Corollary 7]{BoRi_short_Kesten} and Tassion \cite{Tassion}. 

 RSW immediately implies a bound on a `one-ar' event: 
\begin{theorem}
\label{thm: RSW one-arm}
For the Bernoulli percolation with $p=1/2$ on the square lattice there is $\delta>0$ such that for all $s<t$ the probability that there is a percolation cluster connecting the boundary of $[-s,s]^2$ to the boundary of $[-t,t]$ is bounded by $(s/t)^\delta$.
\end{theorem}

\begin{remark}
Theorem \ref{thm: RSW one-arm} immediately implies that a.s. there are no infinite clusters. Hence RSW type inequalities can be used to identify the critical probability $p_c$. One should expect that for $p<p_c$ the crossing probability goes to zero as the scale goes to infinity. For $p>p_c$ we expect it to go to $1$. We will discuss this in more detail in the next section.
\end{remark}

The first RSW result for nodal domains of Gaussian fields was obtained by Beffara and Gayet \cite{BeGa}. They have combined Tassion's very general approach to RSW with quasi-independence estimates that we have discussed in the previous section to obtain the following RSW-type result:
\begin{theorem}[Beffara-Gayet  Theorem 1.1 \cite{BeGa}]
\label{thm: RSW BF}
Let $f$ be a stationary centred Gaussian field in $\R^2$ with the covariance kernel $K$. We assume that $f\in C^4$ a.s., $K(x_1,x_2)=K(x_1,-x_2)$, $K(x)\ge 0$ and $K(x)=o(|x|^{-325})$ as $|x|\to\infty$. Then for any quad (i.e. a smooth bounded open domain $\Omega\subset \R^2$ with  two marked disjoint boundary arcs $\gamma_1$ and $\gamma_2$) there is $c>0$ such that for all $R>0$
\[
\P\brb{\{x\in R\cdot\Omega: f(x)>0 \} \text{ connects } R\cdot\gamma_1 \text{ and } R\cdot\gamma_2}\ge c.
\]
If $R$ is sufficiently large, then the same is true for $\{f=0\}$.
\end{theorem}
\begin{remark}
 Their main aim was to prove RSW for the Bargmann-Fock field, which clearly satisfies all these conditions. Since its covariance decays faster than any polynomial, there was no need to optimize the covariance decay rate.
\end{remark}

The main contribution to the exponent $325$ is due to strong assumptions in their quasi-independence argument. Better quasi-independence results give better assumptions on the covariance. In \cite[Theorem]{BeMu_discr} the exponent was reduced to $16$ and a bit later it was reduced to $4$ in \cite[Theorem 1.1]{RiVa}. It seems that $4$ is the best exponent that could be obtained using quasi-independence methods. Later, Muirhead and Vanneuville \cite{MuVa} reduced the exponent to $2$ using very different ideas. We will discuss it in the next section.

Although \cite{BeKe} gives strong numerical evidence that there is RSW for the random plane wave, current methods are not sufficient to tackle this case. The covariance kernel decays like $|x|^{-1/2}$ which is too slow for quasi-independence results. Moreover, it is not positive, so there is no FKG inequality (see Theorem \ref{thm Pitt}) which is heavily used in the proof of RSW.

A bit later the argument of Beffara and Gayet was extended to the Kostlan ensemble. 
\begin{theorem}[Beliaev-Muirhead-Wigman \cite{BMW}]
\label{thm: RSW Kostlan}
Let $f_n$ be the Kostlan ensemble of degree $n$. For any quad on the unit sphere that does not contain antipodal points, there is $c>0$ such that the probability that $\{f_n>0\}$ crosses it is at least $c$ for all $n$. The same is true for $\{f=0\}$ is $n$ is sufficiently large. 
\end{theorem}
Moreover, the constant $c$ is uniform in the quad, provided that we can control its `shape'. The precise statement can be found in \cite[Definition 1.4]{BMW}.

We would like to point out a few important differences between Theorems  \ref{thm: RSW BF} and \ref{thm: RSW Kostlan}. The second theorem (i) deals with percolation on the sphere which is a bit more complicated than the planar one (ii) the estimate is uniform over a class of different fields (iii) the covariance kernel of the Kostlan ensemble can be negative for odd $n$ (iv) Tassion-style argument is applied not to a discretized field but directly to the continuous field.

\begin{remark}
As we have mentioned in the Remark \ref{rem: Kostlan projective variety} the Kostlan ensemble is a natural model for a `typical' homogeneous polynomial of degree $n$ and its nodal lines are `typical' real algebraic varieties. Theorem \ref{thm: RSW Kostlan} gives strong support to the percolation conjecture for the Kostlan ensemble. This allows us to conjecture that a `typical' projective variety converges to the Conformal Loop Ensemble CLE($6$) as $n$ tends to infinity. This gives an interesting connection between algebraic geometry and Conformal Field Theory. 
\end{remark}

In \cite{BeGa_no_FKG} Beffara and Gayet used perturbation techniques to obtain RSW for a class of models without FKG. They have argued that any model which in a certain sense close to a model with RSW has RSW as well. In particular, this can be applied to Gaussian fields that decorrelate quickly. In a sense this is close to an idea used in \cite{BMW}: if negative correlations are very small, then one can still obtain good low bounds on the intersection of crossing events.

\subsection{Phase transition}

One can think that a field $f$ defines a landscape and the level $\ell$ defines a water level. When $\ell$ is a very large negative number, we should expect that there is a huge `continent' with small `lakes'. As the water level is rising, these lakes are getting larger and larger. At $\ell=0$ there is a balance between the land and the water and after that we switch to an `ocean with islands' regime. When $\ell\to\infty$, these islands are getting smaller and eventually disappear. A similar picture happens in percolation when $p$ is increasing from $0$ to $1$. The phenomenon of switching from `land with lakes' to `ocean with islands' is called a phase transition.

\begin{figure}
\includegraphics[width=0.3\textwidth]{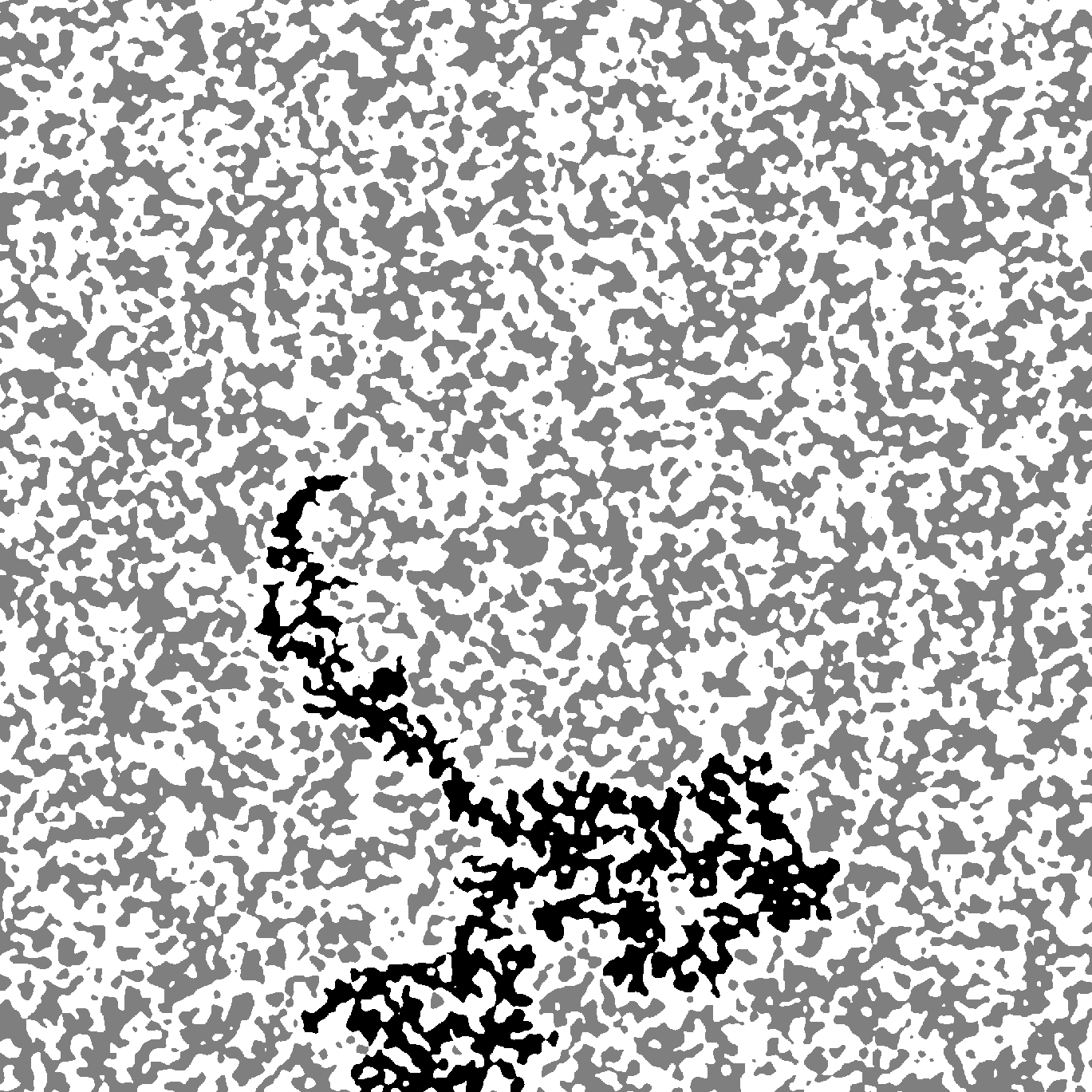}
\hspace{0.03\textwidth}
\includegraphics[width=0.3\textwidth]{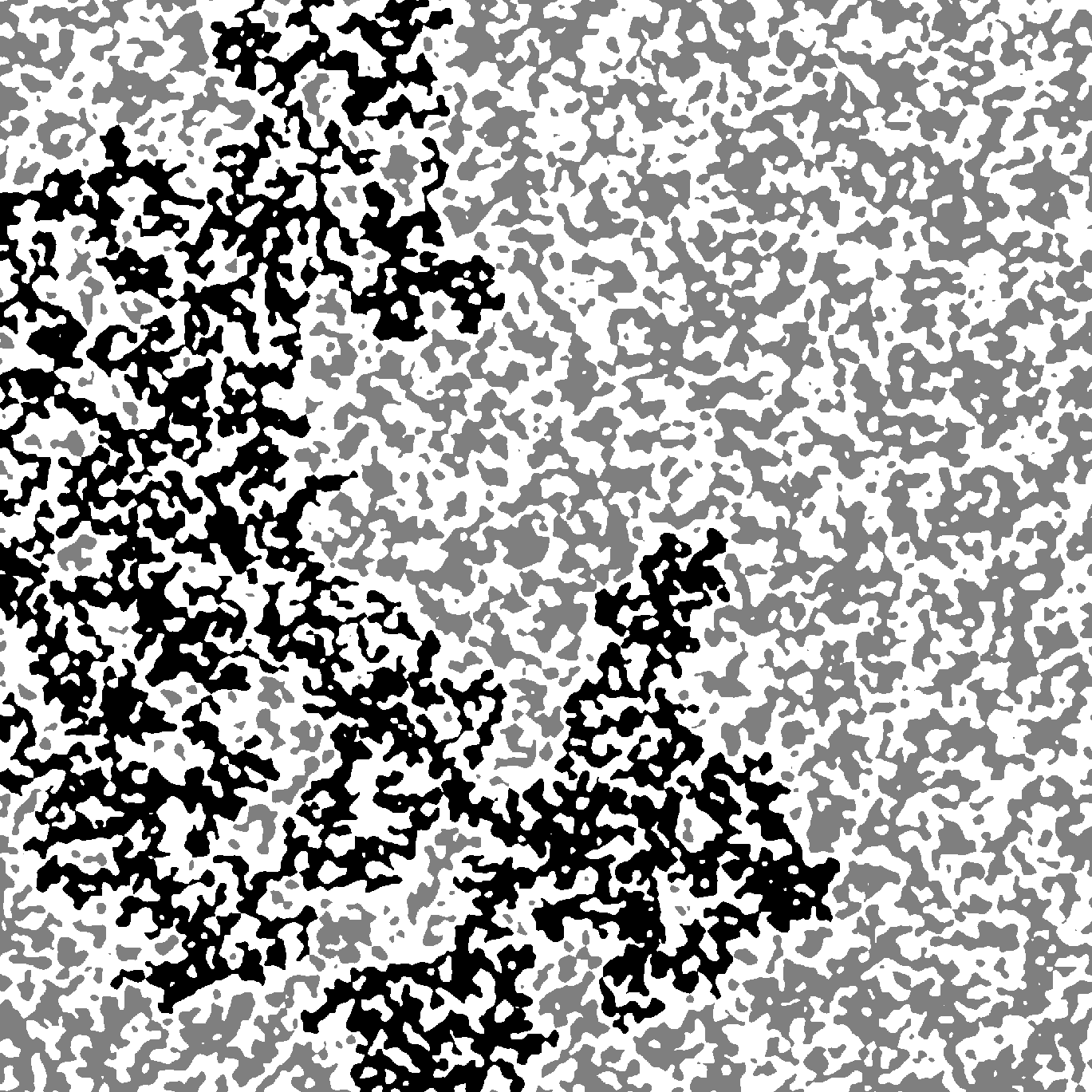}
\hspace{0.03\textwidth}
\includegraphics[width=0.3\textwidth]{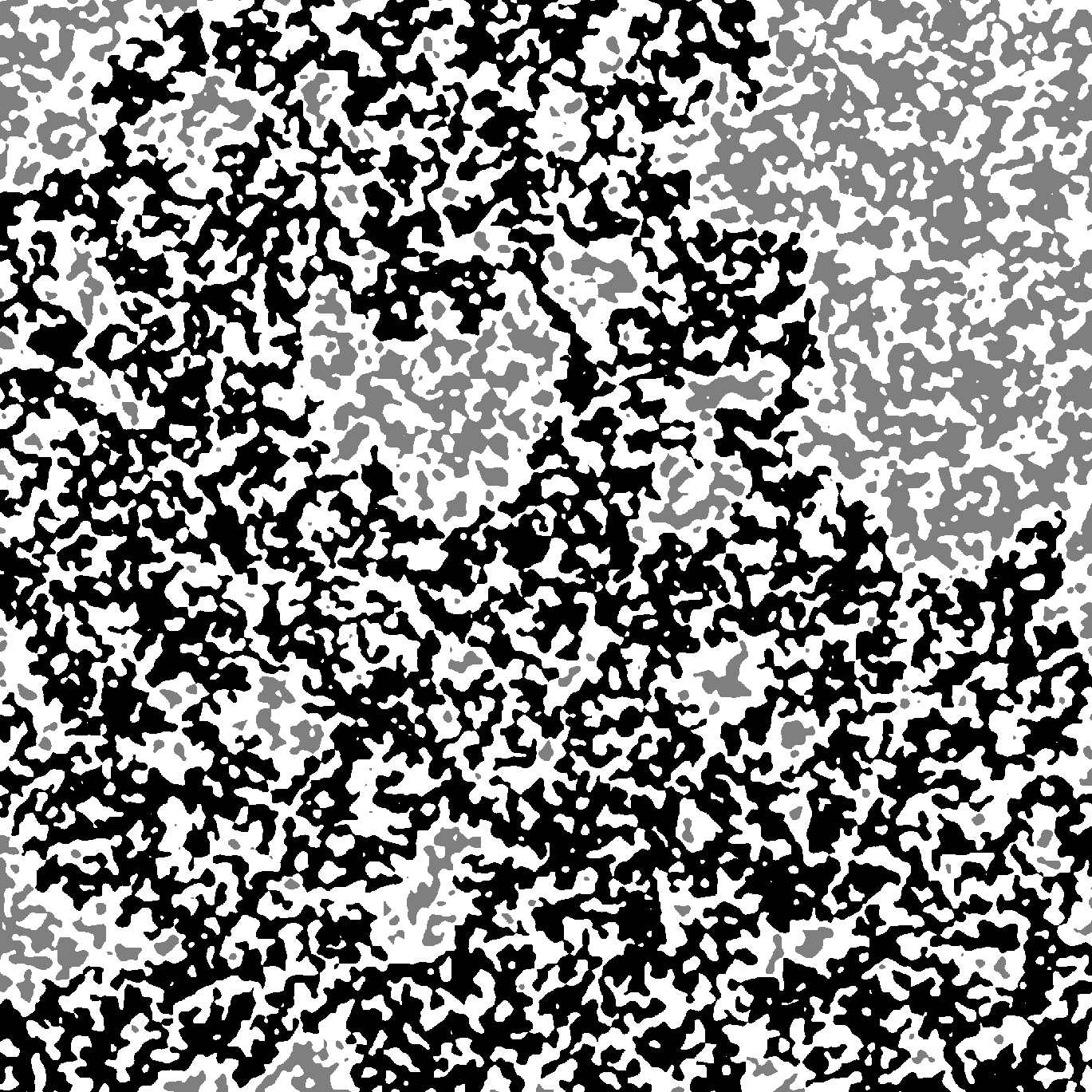}
\caption{Three different level sets of the same sample of the Bargmann-Fock field. In all figures grey is the set above the level, white is below the level and black is the largest component of the excursion set. The levels are $\ell=0.1$ (left), $\ell=0$ (center) and $\ell=-0.1$ (right). As the level is decreasing, the largest component changes from small, to thin but comparable to the size of the domain, to a dense one.} 
\label{fig: BF phase transition}
\end{figure}

It should not be a surprise that the first result of this type has been proved for the Bargmann-Fock field.  See Figure \ref{fig: BF phase transition} for an example of a phase transition.

\begin{theorem}[Rivera-Vanneuville \cite{RiVa_BF_phase_transition}]
\label{thm: RV phase transition}
Let $f$ be the Bargmann-Fock field, then a.s the set $\{f\ge \ell\}$ has a unique unbounded component when $\ell<0$ and a.s. has no unbounded components if $\ell\ge 0$. Moreover, when $\ell<0$ the crossing probability at scale $R$ (i.e. the probability that there is a curve connecting left and right sides of $[0,\rho R]\times [0,R]$ such that $f>\ell$ along this curve) is at least $1-\exp(-cR)$, where $c>0$ depends of $\rho$ and $\ell$ but not on $R$.
\end{theorem} 
We would like to make a few remarks:
\begin{itemize}
\item We use slightly different notations, Rivera and Vanneuville considered sets $\{f\ge -\ell\}$ that are increasing in $\ell$, whereas our sets $\{f \ge \ell\}$ are decreasing. 
\item This result was stated for the Bargmann-Fock field, but in fact the proof of the first part of the theorem is also valid for `nice' fields with sufficiently fast decay of correlations (faster than $|x|^{-5}$). See \cite[Theorem 2.11]{RiVa} for the complete list of assumptions.  The last part of Theorem \ref{thm: RV phase transition} requires super-exponential decay of correlations.
\item The fact that all nodal domains are bounded follows from RSW. By monotonicity the same is true for all $\ell>0$. 
\item The first part of the theorem is similar to Kesten's theorem about phase transition for Bernoulli percolation on $\Z^2$ (see Theorem \ref{thm: Kesten}). The second part is similar to Kesten's exponential decay, see \cite[Theorem 2]{Kesten}. 
\end{itemize}

This paper uses all the standard techniques and assumptions that we have discussed in the previous section. They are needed to have RSW, to control box crossings and to have FKG. Beyond this, Rivera and Vanneuville used one more tool that originates in theoretical computer science: analysis of influences.

To explain the main ideas we first have to introduce some notations. Let $\P_p^n$ be the product $p$-Bernoulli measure on $X=\{0,1\}^n$ (clearly the Bernoulli percolation on a graph with $n$ bonds can be identified with such measure). Let $A$ be an event, i.e. $A\subset  X$. The influence of $i$-th coordinate $\mathrm{Infl}_i^p(A)$ is the probability that changing the value of $i$-th coordinate will change $\id_A$. In other terms, this is the probability that the $i$-th coordinate is \emph{pivotal} for the event $A$.  The first important ingredient is the Russo differential formula
\begin{equation}
\label{eq: Russo formula}
\frac{d}{dp}\P^n_p\brb{A}=\sum_{i=1}^n \mathrm{Infl}_i^p(A),
\end{equation}
where $A$ is any increasing event (typical example: any crossing event). The right hand side can be also interpreted as the expected number of pivotal coordinates. This formula was obtained independently by Russo in \cite[Lemma 3]{Russo} and Margulis \cite[Lemma 2.6]{Margulis}. See also \cite[Section 2.4]{Grimmett} for a more detailed discussion. It has been noted even in the early work of Margulis that this result could be used to show sharp transition. 

The second ingredient is the Kahn–Kalai–Linial (KKL) theorem \cite{KKL,BKKKL} that gives a lower bound on the sum of influences in terms of the maximal influence. Roughly speaking, it tells us that if all influences are small, then their sum is large. This approach to proving Kesten's theorem was used in \cite{BoRi_short_Kesten}. To be more precise, they have used a sharp threshold result of Friedgut and Kalai \cite{FrKa} based on KKL. See also \cite{BoRi_Voronoi_critical} where these ideas are used in the more relevant context of Voronoi percolation. 

For more information about the harmonic analysis of Boolean functions, noise sensitivity and other related results we refer to very good lecture notes by Garban and Steif \cite{GaSt}. 

One of the main contributions of River and Vanneuville is the adaptation of these ideas to the field setting (in particular, they extend the Russo formula and KKL from the Bernoulli product measure to  Gaussian measures. 

A bit later Muirhead and Vanneuville \cite{MuVa} extended this result to a much larger class of fields. They work in the setting of fields having spectral density, hence representable as a convolution of the white noise with a certain kernel $q$ satisfying $K=q*q$. They obtained the following result:

\begin{theorem}[Theorems 1.6 and 1.11 \cite{MuVa}]
Let $f=W*q$ be a stationary Gaussian field satisfying natural smoothness, symmetry and regularity assumptions. We also assume that $q\ge 0$ and that $q$ and its derivatives of order up to $3$ decay faster than $|x|^{-\beta}$ for some $\beta>2$, then there is sharp phase transition for the excursion sets $\{f\ge \ell\}$. Namely, (i) for $\ell>0$ all components are bounded and crossing probabilities on scale $R$ decays exponentially fast as $R\to\infty$ (ii) there is an RSW estimate for $\ell=0$ (in particular, all components are bounded) (iii) for $\ell<0$ there is a unique unbounded component and crossing probabilities are exponentially close to $1$.
\end{theorem} 
Note that the decay rates of $K$ and $q$ are essentially the same, so the second part gives RSW under a weaker decay assumption than \cite[Theorem 1.1]{RiVa}. 

There are few important improvements compared to all previous results. First, is that they do not use discretization of the field, instead, similar to \cite{BMR}, they approximate the field by a smooth field such that the corresponding Hilbert space is finite-dimensional. Secondly, instead of quasi-independence, they use a `sprinkled' version. Namely, for crossing events $A$ and $B$ they compare $\P[f\in A]\P[f\in B]$  to $\P[f-\epsilon\in A\cap B]$ instead of $\P[f\in A\cap B]$. Sprinkling was also used in \cite{RiVa_BF_phase_transition}. For use of sprinkling in percolation we refer to the discussion in Section 2.6 of \cite{Grimmett} and \cite[Theorem 2.46]{Grimmett}. See also \cite[Lemma 4.2]{ACCFR}. Finally, they adapted some new ideas that were recently used by Duminil-Copin,  Raoufi and  Tassion to study lattice models and different versions of percolation \cite{DCRT19a,DCRT19b,DCRT20}.

To be more precise, they have used a version of OSSS  inequality \cite{OSSS} that gives a low bound on the sum of influences. Without going into technicalities, it can give a bound on the sum of influences in terms of \emph{revealments}. 

Let $A\subset X$, where $X=E^n$ is a product space equipped with a probability product measure.   We consider a \emph{random} algorithm $\mathcal{A}$.  The algorithm reveals coordinates of $\omega$ until it is able to determine whether $\omega \in A$ or not. The revealment of $i$-th coordinate is the probability that this coordinate will be revealed by $\mathcal{A}$. Essentially OSSS inequality tells us that if all revealments are small, then the sum of influences is large, hence the derivative in the Russo formula \eqref{eq: Russo formula} is large.  

\begin{figure}
\includegraphics[width=0.8\textwidth]{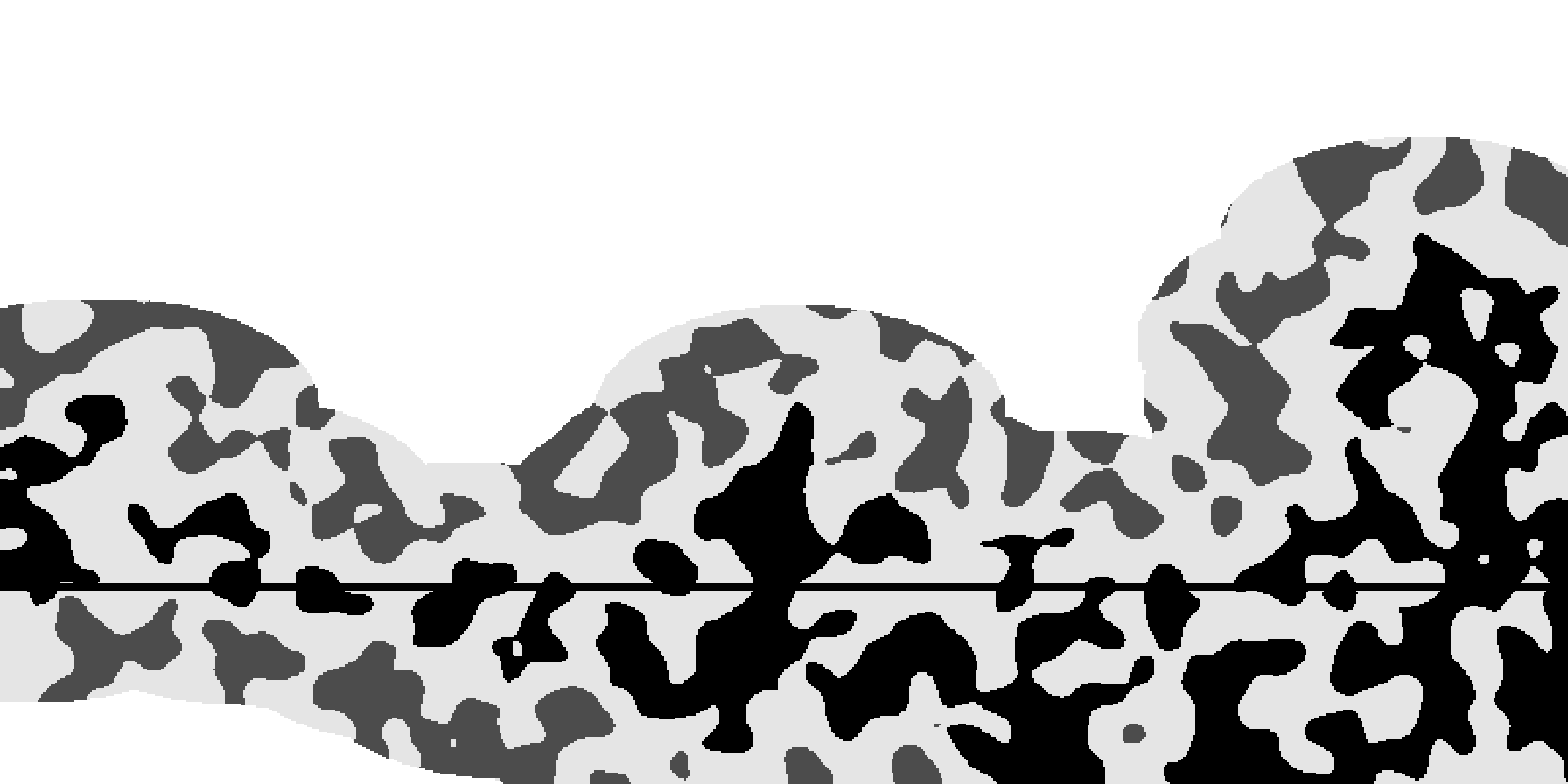}
\caption{An example of the revealment algorithm. Light grey is where $f>0$, dark grey is where $f<0$. All negative components intersecting random horizontal lines are black. Only the part of the field revealed by the algorithm is shown.  To be more precise, to reveal this part of the field we have to reveal all white noise coordinates within a finite distance.}
\label{fig: revealment}
\end{figure}

Muirhead and Vanneuville approximate the field by a convolution of the truncated kernel $q$ with a discrete white noise. The random algorithm works approximately like this: we start with the white noise in a rectangle. We first choose a horizontal line and reveal all white noise coordinates around it so that the values of $f$ on the line are known (there are finitely many of them since $q$ is compactly supported). After that, we explore all negative domains intersecting the line. If one of them crosses vertically, then there is no positive horizontal crossing, otherwise there is a positive crossing. See Figure \ref{fig: revealment} for an example of what is revealed by this algorithm. A white noise coordinate is revealed if it is close to one of these negative domains. This probability is controlled by the one-arm exponent, so RSW gives an upper bound on revealments.  

The assumption in the main result of \cite{MuVa} were slightly improved by Rivera \cite{Rivera_Talagrand_inequality}.  In particular, the strong positivity assumption $q\ge 0$ which is not very natural and hard to verify was replaced by a more natural weak positivity assumption $K\ge 0$. What is probably more important, the argument of Rivera uses a Gaussian Talagrand inequality. This is arguably a much more natural tool in the case of continuous fields than KKL or OSSS inequalities used before. 

The assumptions needed for the phase transition were weakened further in \cite[Theorem 1.3]{MRV}.  Probably the most amazing thing is that their result uses extremely weak assumptions on the covariance decay rate. They only assume that $K$ decays faster than $(\log\log(|x|))^{-\beta}$ for some $\beta>2$. In particular, this result is valid for the random plane wave where the covariance kernel decay so slowly that we are unable to establish quasi-independence and the kernel is not positive, so there is no FKG and negative correlations are not negligible.

To deal with some of these issues the authors have modified Tassion's argument \cite{Tassion} (see also \cite{K-ST,BeGa, BMW}) so that it uses sprinkling to glue together crossing events instead of FKG. The resulting argument is much more robust and uses weaker assumptions but it does not work at all at criticality. So \cite{MRV} establishes the sharp phase transition but not RSW at criticality. 

The second ingredient is the use of Chatterjee's approach to supercontractivity. This is somewhat similar to the use of Talagrand's inequality in \cite{Rivera_Talagrand_inequality}.  The authors used it to show the following criterion for the sharp threshold. For an increasing event $A$ one can define the threshold height $T_A$, which is the first level at which the event happens. For example, let $A$ be the crossing event. Clearly, for very large $\ell$ the excursion set $\{f>\ell\}$ does not cross and for a very large negative $\ell$ it does cross. By monotonicity, there is the threshold $\ell=T_A$ at which the excursion set switches between crossing and non-crossing. At the level $T_A$ there is a pivotal place where the change in the topology happens. This place is called the \emph{threshold location}. The criterion states (in a quantified way) that the sharp threshold is equivalent to the delocalization of the threshold location. On a heuristic level, it is similar to the revealment approach. If there is a place that where the threshold is localized with high probability, then this place should be pivotal for crossing and any algorithm discovering the event $A$ should reveal it with high probability. 

Severo extended these results to higher dimensions. Theorem 1.2 \cite{Severo} states that there is a sharp phase transition for a stationary field in $\R^n$ satisfying the usual regularity, symmetry and smoothness assumptions provided that its covariance kernel is strongly positive ($q\ge 0$) and decays faster than $|x|^{\beta}$ for some $\beta>n$. The argument is based on the comparison between the field and its approximation by a finite-range field. 

There are several recent results that are not directly related to the percolation conjecture, but they also deal with phase transitions and connectivity probabilities in sub- and supercritical regimes in Gaussian fields. They consider these questions for a large class of \emph{strongly correlated} Gaussian fields in $\Z^n$ and $\R^n$. A typical example would be the discrete Gaussian Free Field in dimensions $n\ge 3$ where the covariance kernel decays like $|x|^{d-2}$. It was also shown that there is a RSW-type estimate for very strongly but positively correlated fields. A typical example would be a field with the Cauchy covariance kernel $1/(1+|x|^2)^{\alpha/2}$ with $\alpha<2$. This gives some indication that there are scaling limits beyond the Bernoulli percolation universality class. For more information we refer to \cite{DGRS,GRS,MuSe,Muirhead22} and references wherein. 

\section{What we do not know}

Let us start with a technical question that we still find interesting and important. It could be thought of as questions about reproducing kernel Hilbert spaces (RKHS). As we have seen before, quite often it is useful to have quantified finite-dimensional approximations. This leads to the following question
\begin{question}
Let $H$ be an RKHS made of $C^k$ functions. What is the best finite-dimensional approximation of $H$? That is, we would like to find an $n$-dimensional Hilbert space $H_n$ such that functions from $H$ could be well approximated in $C^k(B(R))$ norm by functions from $H_n$. One can try to optimize the dimension $n$ and how good is the approximation in terms of the radius $R$. 
\end{question}

There is also a similar question about the localization of orthonormal bases of an RKHS $H$. Recall that the corresponding Gaussian field (isonormal process) can be written as $\sum a_k \phi_k$ where $\{\phi_k\}$ is an orthonormal basis and $a_k$ are independent Gaussian random variables. Its covariance kernel is the reproducing kernel $K$ of $H$. Truncating the series after $n$ terms we obtain a field $f_n$ with a covariance kernel $K_n(x,y)=\sum_{k=1}^n \phi_k(x)\phi_k(y)$. Clearly, $K_n\to K$ as $n\to\infty$. We would like to understand the rate of convergence. 
\begin{question}
Given the ball or radius $R$ and $\epsilon>0$, what is the best $n$ such that $f_n$ is $\epsilon$-close to $f$ in $C^k(B(R))$ norm? Equivalently, what is the best $n$ such that all basis functions except $n$ are very small in $B(R)$? 
\end{question}

Here, we can look at our two main examples: the random plane wave and the Bargmenn-Fock field. For the random plane wave we have series \eqref{eq: RPW series}. Since Bessel functions $J_n$ decay exponentially fast in $n$ it can be shown that truncating the series at $cR$ terms gives an exponentially good approximation in the disc of radius $R$. So here the number of terms grows linearly in $R$. 
Similarly, for the Bargmann-Fock field, we can remove all terms where $n,m>cR$, but the number of such terms is of order $R^2$. So the random plane wave admits much better finite approximations. One can think that this is a manifestation of the rigidity of the RPW. One consequence of this better approximation is that the variance of the number of excursion sets (for a generic non-zero level) in a ball of radius $R$ grows like\footnote{Only the low bond on the growth rate is known rigorously but it is believed to be sharp} $R^3$ \cite{BMM_variance}  for the random plane wave and like $R^2$ for the Bargmann-Fock field \cite{BMM_variance,BMM_CLT}. So this difference is not superficial, it tells us something important about fields.
\begin{question}
Is there a way to see this difference by considering the corresponding covariance kernels without using an explicit orthonormal basis? 
\end{question}

Now, let us return to percolation. 
Despite all the progress in the last five years, there are many questions that are still open, in particular, it is still not quite clear under what assumptions the percolation conjecture should hold. It is also not clear whether a Harris criterion heuristics could be made rigorous. Here we will list several open problems/conjectures where we think that the answer is positive.

\begin{conjecture}[Bargman-Fock percolation]
\label{conj: BF percolation}
For the Bargmann-Fock field in $\R^2$ 
\begin{itemize}
\item Show the existence of one of the universal exponents and that they are equal to the percolation ones.
\item Prove the existence of a conformally invariant scaling limit of any observable. For example, prove Cardy's formula for crossing probability.
\item Show that nodal lines converge to the conformal loop ensemble CLE($6$).
\end{itemize}
\end{conjecture}

Although there are very strong reasons to believe that this conjecture is true, we probably have to understand the Bernoulli percolation better before tackling it. After all, so far these results are only known for the percolation on the triangular lattice. Their proofs are based on Cardy's formula and its proof relies very heavily on the specifics of the triangular lattice. Available proofs are rather simple and elegant (especially the one in \cite{KhSm}) but they are very hard to generalize and extend. 

\begin{conjecture}[Random Plane Wave percolation]
For the RPW in $\R^2$
\begin{itemize}
\item Show that there is no percolation at $\ell=0$ i.e. that all nodal lines are bounded. 
\item Show that the variance of the number of nodal domains in a box of size $R$ grows as $R^2$ (possibly up to a polylogarithmic factor). Note that for the non-zero level it is not true. 
\item Prove some form of quasi-independence of crossing events for nodal lines.
\item Prove RSW.
\item Prove the same statements as in Conjecture \ref{conj: BF percolation}.
\end{itemize}
\end{conjecture}

This conjecture is a bit more debatable, but we believe that there is good numerical support for the conjecture. Proving it will obviously require even more new ideas since even fewer percolation ideas can be adapted to the random plane wave setting since its covariance kernel is slowly decaying and oscillating. 

Beyond these particular forms of the percolation conjecture, there are questions about natural extensions:

\begin{question}
It is easy to believe that the Gaussianity of the fields is not paramount. Clearly, it is a very strong property that allows us to use many tools, but it is easy to imagine that some percolation properties are more universal. For example, if we replace independent normal coefficients in the definition of the Bargmann-Fock field with some other independent random variable, then the series will converge provided the tails are not heavy. Will the nodal domains of such fields behave like the critical percolation? 
\end{question}

The other related class of fields are shot noise fields. Let $\mathcal{P}$ be a Poisson point process. We define a random measure $\mu=\sum a_i \delta_{y_i}$ where the sum is over all points $y_i$ in $\mathcal{P}$ and the coefficients $a_i$ are i.i.d random variables. For an integrable function $g$ we define the shot noise field as the convolution $\mu*g$. It is known \cite{LaMu} that under some assumptions there is a phase transition at $0$ for the excursion sets of such fields. 

\begin{question}
Are nodal domains of shot noise fields behave like the critical percolation?
\end{question}

\begin{acks}[Acknowledgments]
I would like to thank M.~Lifshits for encouraging me to write this survey and S.~Muirhead for valuable comments on the first draft. 

\end{acks}


\bibliographystyle{imsart-number} 
\bibliography{References}       


\end{document}